\documentclass[12pt,draft]{amsart}
\usepackage[mathscr]{eucal}\catcode`\@=11
\usepackage{amssymb}

\advance\headheight by 3pt

\renewcommand{\subjclass}[1]{\thanks{\emph{2000 Mathematics Subject Classification:}~#1}}
\renewcommand{\keywords}[1]{\thanks{\emph{Keywords and Phrases:}~#1}}

\theoremstyle{plain}
\newtheorem{theorem}{Theorem}[section]
\newtheorem{lemma}{Lemma}[section]
\newtheorem{corollary}{Corollary}[section]
\newtheorem{proposition}[lemma]{Proposition}

\def\teto#1{\setbox\z@\hbox{${#1\vphantom k}$}\hbox{%
 \hbox{\lower2\ex@\hbox{\lower\dp\z@\hbox{\vbox{\hrule
 \hbox{\vrule\hskip2\ex@\vbox{\vskip2\ex@\box\z@\vskip1\ex@}%
 \hskip2\ex@\vrule}}}}}}}
\catcode`\@=\active

\newcommand{\Q}{\mathbb{Q}}
\newcommand{\Z}{\mathbb{Z}}
\newcommand{\NN}{\mathbb{N}}
\newcommand{\K}{\Bbbk}

\newcommand{\OO}{\mathcal{O}}
\newcommand{\M}{\mathcal{M}}

\newcommand{\F}{\mathcal{F}}

\newcommand{\CC}{\mathcal{C}^*}
\newcommand{\CCC}{\mathcal{C}}

\newcommand{\PP}{\mathbb{P}}

\newcommand{\vp}{\varphi}
\newcommand{\ve}{\varepsilon}
\newcommand{\al}{\alpha}
\newcommand{\be}{\beta}
\newcommand{\ga}{\gamma}
\newcommand{\de}{\delta}

\newcommand{\fp}{\frak{p}}
\newcommand{\fa}{\frak{a}}
\newcommand{\fb}{\frak{b}}
\newcommand{\fc}{\frak{c}}
\newcommand{\fd}{\frak{d}}
\newcommand{\A}{\frak{a}}
\newcommand{\fA}{\frak{A}}
\newcommand{\fD}{\frak{D}}
\newcommand{\fP}{\frak{P}}

\newcommand{\til}[1]{\tilde{#1}}

\newcommand{\ketsor}[2]{\genfrac{}{}{0pt}{2}{#1}{#2}}

\newcommand{\ordv}{{\rm ord}_v}

\DeclareMathOperator{\Gal}{Gal}

\DeclareMathOperator{\Tr}{Tr}

\newcommand{\simm}[1]{\overset{#1}{\sim}}
\newcommand{\simmm}[1]{\overset{#1}{\approx}}
\newcommand{\kacs}[1]{\overset{#1}{\prec}}

\newcommand{\kdots}{,\ldots ,}
\renewcommand{\eqref}[1]{(\ref{#1})}
\begin{document}

\title[Binary forms of given degree and given discriminant]
{On the number of equivalence classes of binary forms of given degree
and given discriminant}


\author[A. B\'erczes]{Attila B\'erczes}
\thanks{The research was supported in part by the Hungarian Academy of
Sciences (A.B.,K.G.), the Netherlands Organization for Scientific
Research (A.B.,J.-H.E.,K.G.), by
grants F34981 (A.B), N34001 (A.B.,J.-H.E.,K.G.) T42985 (A.B., K.G.) and
T38225 (A.B., K.G.)  of the Hungarian National Foundation for Scientific
Research and by the FKFP grant 3272-13066/201 (A.B.)}
\address{A. B\'erczes \newline
         \indent Institute of Mathematics, University of Debrecen \newline
         \indent Number Theory Research Group, Hungarian Academy of Sciences and \newline
         \indent University of Debrecen \newline
         \indent H-4010 Debrecen, P.O. Box 12, Hungary}
\email{berczesa\char'100math.klte.hu}

\author[J.-H. Evertse]{Jan-Hendrik Evertse}
\address{J.-H. Evertse \newline
         \indent Mathematical Institute, Universiteit Leiden \newline
         \indent P.O. Box 9512, NL-2300 RA Leiden, The Netherlands}
\email{evertse\char'100math.leidenuniv.nl}

\author[K. Gy\H{o}ry]{K\'{a}lm\'{a}n Gy\H{o}ry}
\address{K. Gy\H{o}ry \newline
         \indent Institute of Mathematics, University of Debrecen \newline
         \indent Number Theory Research Group, Hungarian Academy of Sciences and \newline
         \indent University of Debrecen \newline
         \indent H-4010 Debrecen, P.O. Box 12, Hungary}
\email{gyory\char'100math.klte.hu}

\maketitle

\centerline{\textit{To Professor Robert Tijdeman on his 60th birthday}}

\section{Introduction}

In the present paper we give explicit upper bounds for the number of
equivalence classes of binary forms of given degree and discriminant,
and for the number of equivalence classes of irreducible
binary forms with given invariant order.

Two binary forms $F,G\in \Z [X,Y]$ are called equivalent if there
is a matrix ${a~b\choose c~d}\in {\rm GL}_2(\Z )$
such that $G(X,Y)=F(aX+bY,cX+dY)$. 
Denote by $D(F)$ the discriminant of a binary form $F$,
and by $\OO_F$ the
invariant order of an irreducible binary form $F$.
We recall the definition of the invariant order of $F$ which is less familiar.
Write $F(X,Y)=a_0X^r+a_1X^{r-1}Y+\cdots +a_rY^r$
and let $\theta_F$ be a zero of $F(X,1)$. Then $\OO_F$ is defined
to be the $\Z$-module with basis $1$, $a_0\theta_F$,
$a_0\theta_F^2+a_1\theta_F$,
$a_0\theta_F^3+a_1\theta_F^2+a_2\theta_F$,$\ldots$,
$a_0\theta_F^{r-1}+a_1\theta_F^{r-2}+\cdots +a_{r-2}\theta_F$;
this is indeed an order, i.e., closed under multiplication.
It is well-known that two equivalent binary forms have the same
discriminant. Further, two equivalent irreducible binary forms
have the same invariant order. The discriminant $D(\OO_F)$
of $\OO_F$ is equal to $D(F)$ (see 
\cite{Naka}, \cite{Simon} for a verification
of these facts).
Consequently, if $K=\Q (\theta_F)$,
then $D(F)=c^2D_K$, where $D_K$ is the discriminant of $K$ and
$c=[\OO_K:\OO_F]$ is the index of $\OO_F$ in the ring of integers
$\OO_K$ of $K$.

By classical results of Lagrange, Gauss ($r=2$) and Hermite $(r=3)$,
the binary forms $F\in\Z [X,Y]$ of degree $r \leq 3$
with a given discriminant $D \ne 0$ lie in finitely many
equivalence classes, and these classes can be effectively
determined. This finiteness
theorem was generalized for the case $r\geq 4$ by
Birch and Merriman \cite{BirMer} in an ineffective form, and
later by Evertse and Gy\H ory \cite{EGy5} in an effective form.
Moreover, the theorem remains true without fixing the degree $r$;
see \cite{Gy14}.
An immediate consequence is that if $\OO$
is a given order of some number field, then the irreducible binary forms
$F\in \Z [X,Y]$ with $\OO_F=\OO$ lie in finitely many equivalence
classes. From a result of Delone and Faddeev \cite[Chap.II, \S15]{DelFad}
it follows that for each cubic order $\OO$ there is precisely one
equivalence class of irreducible binary cubic forms $F\in\Z [X,Y]$
such that $\OO_F =\OO$. For degree larger than $3$ this is no longer true:
Simon \cite{Simon} gave examples of number fields $K$ of
degree $4$ and of arbitrarily
large degree whose ring of integers $\OO_K$ can not be represented
as $\OO_F$ for any irreducible binary form $F$.

In the present paper, we prove the following results:
\\[0.15cm]\indent
{\bf 1)} Let $\OO$ be an order whose quotient field has degree $r\geq 4$
over $\Q$. Then the irreducible binary forms $F\in\Z [X,Y]$ with $\OO_F\cong\OO$
lie in at most $2^{24r^3}$ equivalence classes.
\\[0.15cm]\indent
{\bf 2)} Let $K$ be an algebraic number field of degree $r\geq 3$
and let $c$ be a positive integer.
Then for every $\varepsilon >0$
the set of irreducible binary forms
$F\in \Z [X,Y]$ such that $K=\Q (\theta_F)$ for some zero $\theta_F$
of $F(X,1)$ and such that $D(F)=c^2D_K$ is contained in the union of
at most
$\al (r,\varepsilon)c^{\frac{2}{r(r-1)}+\varepsilon}$
equivalence classes; here $\al (r,\varepsilon )$ depends only on $r$
and $\varepsilon$.
We show that in this upper bound the exponent of $c$
cannot be replaced by a quantity smaller than $\frac{2}{r(r-1)}$.
\\[0.15cm]\indent
More generally, we prove analogues of 1) and 2) for binary forms
having their coefficients in the ring of $S$-integers of a number field.
Further, we prove a generalization of 2) for reducible binary forms.
Our precise results are stated in Section 2
(Theorems 2.1, 2.2 and 2.3).
 Our approach is similar to that of Birch
and Merriman \cite{BirMer}, with the necessary modifications.
In our proofs we use among other things an upper bound by Beukers
and Schlickewei \cite[Theorem 1]{BeSchl} for the numbers of solutions of
the equation $x+y=1$ in unknowns $x,y$ from a multiplicative group
of finite rank.
\\[0.5cm]

\section{Statements of the results}\label{S1}

\textbf{Terminology.}
Before stating our results we introduce the necessary terminology.
Let $F(X,Y)=a_0X^r+a_1X^{r-1}Y+ \dots +a_rY^r$ be a binary form.
Writing $F$ as 
\[
F(X,Y)=\lambda\prod_{i=1}^r (\al_iX-\be_i Y)
\]
we may express the discriminant of $F$ as
\begin{equation}\label{1.1}
D(F)=\lambda^{2r-2}\prod_{1\leq i < j \leq r}(\al_i\be_j-\al_j\be_i)^2\, .
\end{equation}
This is independent of the choice of $\lambda$ and of the $\al_i,\be_i$.
It is well-known that $D(F)$ is a homogeneous
polynomial of degree $2r-2$ in $\Z [a_0, \dots , a_r]$.
For a matrix
$U=\bigl( \begin{smallmatrix}a&b\\c&d\end{smallmatrix}\bigr)$
we define $F_U(X,Y):=F(aX+bY,cX+dY)$. Then \eqref{1.1} gives
\begin{equation}\label{1.3}
D(F_U)=(\det U)^{r(r-1)}D(F).
\end{equation}

Now let $R$ be an integral domain with quotient field of characteristic 0.
Two binary forms $F,G \in R[X,Y]$ are called $R$-\textit{equivalent}, notation
$F \simm{R} G$, if $G=F_U$ for some matrix
$U \in \text{GL}_2(R)$, i.e., with $\det U \in R^*$.
(If $R=\Z$ we simply speak about equivalence.)
It is then clear from
\eqref{1.3} that for any two binary forms $F,G \in R[X,Y]$ we have
\begin{equation}\label{1.4}
G \simm{R} F \ \ \Rightarrow \ \ D(G)=\ve D(F) \ \ \text{for some} \ \ \ve\in R^*.
\end{equation}

An important invariant of an irreducible binary form $F\in R[X,Y]$ is its
\textit{invariant ring} or \textit{invariant order} $\OO_{F,R}$
(see Simon \cite{Simon}). By an $R$-order of degree $r$
(or just an order of degree $r$ if $R=\Z$) we mean an
integral domain $\OO$ such that $\OO$ is an overring of $R$,
the domain $\OO$ is
finitely generated as an  $R$-module, and the quotient field of $\OO$
has degree $r$ over the quotient field of $R$.

The order $\OO_{F,R}$ (or just $\OO_F$ if $R=\Z$) is defined as follows.
Let $F=a_0X^r+a_1X^{r-1}Y+ \dots +a_rY^r$ be a binary form in $R[X,Y]$
which is irreducible over the quotient field of $R$.
Let $\theta_F$
be a zero of $F(X,1)$. Then $\OO_{F,R}$ is defined to be the $R$-module
with basis
\begin{equation}\label{1.5}
\begin{split}
\omega_1=1, \ \omega_2=a_0\theta_F, \
&\omega_3=a_0\theta_F^2+a_1\theta_F, \dots ,\\
&\omega_r=a_0\theta_F^{r-1}+a_1\theta_F^{r-2}+ \dots + a_{r-2}\theta_F.
\end{split}
\end{equation}
We recall some facts proved by Simon \cite{Simon} about $\OO_{F,R}$.
First $\OO_{F,R}$ is an $R$-order of degree $r$. Second, if
$G$ is another binary form in $R[X,Y]$ then
\begin{equation}\label{1.5b}
F \simm{R} G  \ \ \Rightarrow \ \ \OO_{F,R} \cong \OO_{G,R} \ \ \ \text{(as $R$-algebras)}.
\end{equation}
Third
\begin{equation}\label{1.5a}
D(\omega_1, \dots ,\omega_r)=D(F).
\end{equation}
Here $D(\omega_1, \dots ,\omega_r)$ denotes the discriminant of
$\omega_1, \dots ,\omega_r$, that is the determinant
$\det(\Tr(\omega_i\omega_j)_{1\leq i,j \leq r})$, where
$\Tr$ denotes the trace map from the quotient field of $\OO_{F,R}$
to that of $R$.
\\[0.2cm]\indent
Our results will be established for binary forms having their
coefficients in the
ring of $S$-integers of a number field. Therefore we recall some notions
about such rings.

Let $\K$ be a number field, and
$\{ | . |_v : v \in M_{\K} \}$ be a maximal set of pairwise inequivalent
absolute values of $\K$. We will refer to $M_{\K}$ as the set of places of $\K$.
Let $S$ be a finite subset of $M_{\K}$ containing
all infinite places of $\K$
(i.e., the places $v$ such that $| . |_v$ is archimedean).
Then the ring of $S$-integers
and its unit group, the group of $S$-units are defined by
$$
\OO_S=\{x \in \K : |x|_v \leq 1 \ \text{for} \ v\not\in S \}, \ \ \
\OO_S^*=\{x \in \K : |x|_v = 1 \ \text{for} \ v\not\in S \},
$$
respectively.

Two ideals $\fa$, $\fb$ of $\OO_S$
are said to belong to the same ideal class of $\OO_S$
if there are non-zero $\lambda ,\mu\in\OO_S$ such that
$\lambda\fa=\mu\fb$.
Denote by $h_m(\OO_S)$ the number of ideal classes $\fA$ of $\OO_S$ such that
$\fA^m$ is the class of principal ideals of $\OO_S$.
For a finite extension $K$ of $\K$,
let $\fd_{K/\K,S}$ denote the relative $S$-discriminant, i.e., the ideal of
$\OO_S$ generated by all discriminants $D_{K/\K}(\omega_1, \dots , \omega_r)$,
where $\omega_1, \dots , \omega_r$ runs through all $\K$-bases of $K$ with
$\omega_1, \dots , \omega_r$ integral over $\OO_S$. The absolute norm
of an ideal $\fa$ of $\OO_S$ is defined by $N_S(\fa):=\#\OO_S/\fa$.

Given an irreducible binary form $F \in \OO_S[X,Y]$ we write $\OO_{F,S}$
for its invariant order $\OO_{F, \OO_S}$.
\\[0.3cm]

\textbf{New results.}
Let $\K$, $\OO_S$ be as above.
From results of Birch and Merriman from 1972 \cite{BirMer} (ineffective)
and Evertse and Gy\H{o}ry from 1991 \cite{EGy5} (effective) it follows that
for given $r \geq 2$ and $D\in\OO_S$ with $D\not= 0$,
the binary forms $F\in\OO_S[X,Y]$
with degree $r$ and with $D(F)\in D\OO_S^*$ lie in finitely many $\OO_S$-equivalence classes.
Together with \eqref{1.5a} this implies that for any given $\OO_S$-order
$\OO$, the binary forms $F\in\OO_S[X,Y]$ which are irreducible over $\K$
and for which $\OO_{F,S}=\OO$ lie in finitely many $\OO_S$-equivalence classes.
From a result of Evertse and Gy\H{o}ry \cite[Thm. 11]{EGy7}
%Crelle 1985%
it can be deduced that for a given $\OO_S$-order $\OO$,
the \emph{monic} binary forms $F\in\OO_S[X,Y]$ (i.e., such that $F(1,0)=1$)
with $\OO_{F,S}=\OO$ lie in at most $c(r)^s$ $\OO_S$-equivalence classes,
where $c(r)$ depends only on $r$ and where
$s=\# S$. Our first result extends this
to non-monic binary forms.

\begin{theorem}\label{T1.1}
Let $S \subset M_{\K}$ be a finite set of cardinality $s$, containing
all infinite places. Let $\OO$ be an $\OO_S$-order
of degree $r \geq 3$. Then there are only finitely many $\OO_S$-equivalence
classes of binary forms $F \in \OO_S[X,Y]$ such that $F$ is irreducible in
$\K[X,Y]$ and
\begin{equation}\label{1.8}
\OO_{F,S} \cong \OO \ \ \ (\text{as $\OO_S$-algebras}).
\end{equation}
The number of these classes is bounded above by
\begin{equation}\label{1.9}
\left\{
\begin{aligned}
&2^{24r^3s} \ \ \ \ \ &&\text{if $r$ is odd},\\
&2^{24r^3s}h_2(\OO_S) \ \ \ \ \ &&\text{if $r$ is even}.
\end{aligned}
\right.
\end{equation}
\end{theorem}

In Section \ref{S8} we show that the factor $h_2(\OO_S)$ is
necessary if $r$ is even.

In the next corollary we state the consequence for $\OO_S=\Z$. Recall that
in this case $\K=\Q$ and $\# S=1$.

\begin{corollary}\label{C1.1}
Let $\OO$ be an order of degree $r \geq 3$.
Then the number of equivalence classes of binary forms $F \in \Z[X,Y]$ such that
$F$ is irreducible in $\Q[X,Y]$ and $\OO_{F} \cong \OO$ is at most
$$
2^{24r^3}.
$$
\end{corollary}

We now state our second result. For an ideal $\fa$ of $\OO_S$, denote by
$\omega_S(\A)$ the number of distinct prime ideals $\fp$ of $\OO_S$
with $\fp \mid \A$
(or the number of $v \not\in S$ such that $|x|_v < 1$ for every $x\in \A$).
Further, for an ideal $\A$ of $\OO_S$ and for $\al \in \NN$, denote by
$\tau_{\al}(\A)$ the number of tuples of ideals $(\fd_1, \dots , \fd_{\al})$
of $\OO_S$ such that their product $\prod_{i=1}^{\al}\fd_i$ divides $\A$.
In the theorems below, the ideal of $\OO_S$ 
generated by $a$ is denoted by $[a]$.

Given a finite extension
$K$ of $\K$, we denote by $\F (\OO_S,K)$ the set of binary forms
$F$ such that $F\in\OO_S[X,Y]$, $F$ is irreducible in $\K [X,Y]$,
and there is $\theta_F$ such that $F(\theta_F ,1)=0$ and $K=\K (\theta_F)$.
By Lemma \ref{L3.1} in Section \ref{S3}, for every $F\in \F (\OO_S,K)$
there is an ideal $\fc$ of $\OO_S$ such that
\begin{equation}\label{1.10}
[D(F)]=\fc^2 \cdot \fd_{K/\K,S}.
\end{equation}

\begin{theorem}\label{T1.2}
Let $S$ be as in Theorem \ref{T1.1}, and let $K$ be an extension of
$\K$ of degree $r \geq 3$.
Then for every non-zero ideal $\fc$ of $\OO_S$, there are at most finitely
many $\OO_S$-equivalence classes of binary forms $F \in \F(\OO_S ,K)$
with \eqref{1.10}.
The number of these classes is at most
\begin{equation}\label{1.11}
2^{24r^3(s+\omega_S(\fc))}\cdot \tau_{\frac{1}{2}r(r-1)}(\fc^2)
\left( \sum_{\fd^{\frac{1}{2}r(r-1)} \mid c} N_S(\fd) \right)\cdot h(r,\OO_S)
\end{equation}
where
\[
h(r,\OO_S)=1\ \ \ \mbox{if $r$ is odd,}\qquad
h(r,\OO_S)=h_2(\OO_S)\ \ \mbox{if $r$ is even.}
\]
Here the sum is taken over all ideals $\fd$ of $\OO_S$ such that
$\fd^{\frac{1}{2}r(r-1)}$ divides $\fc$.
\end{theorem}

We give again the consequence for $\OO_S=\Z$. Given a nonzero integer
$a$, denote by $\omega(a)$ the number of distinct primes dividing $a$, and
for $\al \in \NN$ denote by $\tau_{\al}(a)$ the number of tuples
of positive integers $(d_1, \dots , d_{\al})$
such that $\prod_{i=1}^{\al} d_i$ divides $a$.

\begin{corollary}\label{C1.2}
Let $K$ be a number field of degree $r \geq 3$, and let $c$ be a positive
integer. Then the irreducible binary forms
$F\in\Z [X,Y]$, for which
$\Q (\theta_F)=K$ for some zero $\theta_F$ of $F(X,1)$, and for which
$$
D(F)=c^2D_K
$$
lie in at most
$$
2^{24r^3(1+\omega(c))} \cdot \tau_{\frac{1}{2}r(r-1)}(c^2)
\left( \sum_{d^{\frac{1}{2}r(r-1)} \mid c} d \right)
$$
equivalence classes.
\end{corollary}

Theorem \ref{T1.2} will be deduced from Theorem \ref{T1.1} as follows.
Let $S'$ consist of the places in $S$ and those places $v\not\in S$
such that $|x|_v<1$ for every $x\in\fc$. Then if $F\in\F (\OO_S,K)$ satisfies
\eqref{1.10}, then $D(F)\cdot\OO_{S'}=\fd_{K/\K ,S'}$
and so $\OO_{F,S'}=\OO_{S'}$. Now Theorem \ref{T1.1} yields an upper
bound for the number of $\OO_{S'}$-equivalence classes containing
the binary forms $F\in\F (\OO_S,K)$ with \eqref{1.10} and from the
arguments in Section \ref{S3} one obtains an upper bound for the
number of $\OO_S$-equivalence classes containing the forms lying
in a single $\OO_{S'}$-equivalence class.

We state a generalization of Theorem \ref{T1.2} for reducible forms. Let
$K_0$, $K_1$,$\ldots$,$K_t$ be (not necessarily distinct)
finite extensions of $\K$.
Denote by $\F (\OO_S,K_0\kdots K_t)$ the set of binary forms $F$ with the
following properties: there are binary forms $F_0\kdots F_t$
with $F=\prod_{i=0}^t F_i$, such that
$F_i\in\OO_S[X,Y]$, $F_i$ is irreducible in $\K [X,Y]$,
and there is a $\theta_{F_i}$ such that $F_i(\theta_{F_i})=0$ and
$\K (\theta_{F_i})=K_i$ ($i=0\kdots t$).
By Lemma \ref{L3.1} in Section \ref{S3}, for every binary form
$F\in \F (\OO_S ,K_0\kdots K_t)$
there is an ideal $\fc$ in $\OO_S$ such that
\begin{equation}\label{1.12}
[D(F)]=\fc^2 \fd_{K_0/\K, S} \dots \fd_{K_t/\K, S}.
\end{equation}

\begin{theorem}\label{T1.3}
Let $S$ be as in Theorems \ref{T1.1} and \ref{T1.2}, and let
$K_0, K_1\kdots K_t$
be finite extensions of $\K$. Put $r_i:=[K_i:\K]$ $(i=0, \dots ,t)$ and
$r:=r_0+ \dots +r_t$. Assume that $r_0 \geq 3$.
Then for every non-zero ideal $\fc$ of $\OO_S$
there are at most finitely many $\OO_S$-equivalence classes of
binary forms $F\in\F (\OO_S ,K_0\kdots K_t)$ with \eqref{1.12}.
The number of these classes
is at most
\begin{equation}\label{1.13}
2^{24r^3(s+\omega_S(\fc ))}\cdot \tau_{\frac{1}{2}r(r-1)}(\fc^2)
\left( \sum_{\fd^{\frac{1}{2}r(r-1)} \mid \fc} N_S(\fd) \right)
\cdot h(r_0,\OO_S)
\end{equation}
where
\[
h(r_0,\OO_S)=1\ \ \ \mbox{if $r_0$ is odd,}\qquad
h(r_0,\OO_S)=h_2(\OO_S)\ \ \mbox{if $r_0$ is even.}
\]
\end{theorem}

The consequence of Theorem \ref{T1.3} for $\OO_S=\Z$ is as follows.

\begin{corollary}\label{C1.3}
Let $K_0, \dots , K_t$ be number fields. Put $r_i:=[K_i:\Q]$ $(i=0, \dots ,t)$
and $r:=r_0+\dots + r_t$. Assume that $r_0 \geq 3$. Let $c$ be a positive integer.
Then the binary forms $F$ for which
there are irreducible binary forms $F_0\kdots F_t\in\Z [X,Y]$ with
$F=\prod_{i=0}^t F_i$ such that
$K_i=\Q (\theta_{F_i})$
for some zero $\theta_{F_i}$ of $F_i(X,1)$,
and for which
$$
D(F)=c^2D_{K_0} \dots D_{K_t},
$$
lie in at most
$$
2^{24r^3(1+\omega (c))}\cdot \tau_{\frac{1}{2}r(r-1)}(c^2)
\left( \sum_{d^{\frac{1}{2}r(r-1)} \mid c} d \right)
$$
equivalence classes.
\end{corollary}

Unfortunately, our method of proof of Theorem \ref{T1.3}
requires that we have to impose some
unnatural technical conditions on the
binary forms $F$ under consideration, namely that
they factor into
binary forms $F_i$ with coefficients in $\OO_S$
and that $F_0$ has degree $r_0\geq 3$.
If $\OO_S$ is a principal ideal domain (for instance
when $\K =\Q$), then the first condition is no restriction.
For in that case, if a binary form $F\in\OO_S[X,Y]$ is
reducible over $\K$ its irreducible factors can always be chosen from
$\OO_S[X,Y]$. But the latter is not true if  
$\OO_S$ is not a principal ideal domain.

Allowing these technical conditions, we give a relatively simple proof
of Theorem \ref{T1.3} based on Theorem \ref{T1.2} and on a result on
resultant equations (see Proposition \ref{P7.1} in Section \ref{S7})
which may be of
some independent interest.
It may be possible to remove the technical
conditions from Theorem \ref{T1.3}
at the price of more complications.

Theorem \ref{T1.3} implies that the number of $\OO_S$-equivalence classes
of binary forms $F\in\F (\OO_S, K_0\kdots K_t)$ with \eqref{1.12} is at most
\begin{equation}\label{1.14}
\al(\K, S, r_0, \dots , r_t, \ve)N_S(\fc)^{\frac{2}{r(r-1)}+\ve}
\end{equation}
for every $\ve > 0$, where $\al$ depends only on the parameters between the
parentheses. In Section \ref{S8} we will show that the
bound \eqref{1.14} is almost best possible in terms of $N_S(\fc)$ in
the following sense: for each tuple $(K_0, \dots , K_t)$ of finite
extensions of $\K$, there is a sequence of ideals $\fc$ of $\OO_S$ with
$N_S(\fc )\to\infty$, such that the number of $\OO_S$-equivalence classes
of binary forms $F \in \F(\OO_S,K_0\kdots K_t)$ with \eqref{1.12} is at least
$$
\be N_S(\fc)^{\frac{2}{r(r-1)}},
$$
where $\be$ is a positive constant independent of $\fc$.
\\[0.5cm]

\section{Preliminaries}\label{S2}

\setcounter{equation}{0}

In our proofs it will be necessary to keep track not only
of binary forms but also of their zeros. To facilitate this, we introduce
below so-called augmented forms, 
which are tuples consisting of a binary form
and of some of their zeros.

Given a field $K$, we define $\PP^1(K):=K\cup \{ \infty \}$. Every matrix
$A = \bigl( \begin{smallmatrix} a&b\\c&d\end{smallmatrix} \bigr) \in \text{GL}_2(K)$ induces a
projective transformation
$$
\langle A \rangle : \PP^1(K) \to \PP^1(K) : \xi \mapsto \frac{a\xi+b}{c\xi+d}
$$
(with the usual rules $(a\xi+b)/(c\xi+d)=\infty$ if $c \ne 0$ and $\xi=-d/c$;
$(a\infty+b)/(c\infty+d)=a/c$ if $c \ne 0$ and $\infty$ if $c=0$).
Thus, two
matrices $A,B\in \text{GL}_2(K)$ induce the same projective transformation
if and only if $B=\lambda A$ for some $\lambda \in K^*$. 
\\[0.2cm]

Now let $\K$ be a number field which is fixed henceforth. Let $K$ be a finite
extension of $\K$. An \textit{augmented} $K$-\textit{form} is a pair
$F^*=(F,\theta_F)$
consisting of a binary form $F$ which is irreducible in $\K[X,Y]$,
and $\theta_F \in K$ such that $F(\theta_F,1)=0$ and $\K(\theta_F)=K$.
We agree that $\K(\infty)=\K$ and that for every $c \in \K^*$, $(cY,\infty)$
is an augmented $\K$-form.

Let  $K_0, \dots , K_t$ be a sequence of finite extensions of $\K$. An
\textit{augmented}
$(K_0,\dots , K_t)$-\textit{form} is a tuple $F^*=(F,\theta_{0,F}, \dots ,\theta_{t,F})$ with the property that there are binary forms
$F_0\kdots F_t$,
such that $F=\prod_{i=0}^t F_i$, and
$(F_i, \theta_{i,F})$
is an augmented $K_i$-form for $i=0\kdots t$. 
We define the discriminant and degree of $F^*$ by
$D(F^*):=D(F)$, $\deg F^* := \deg F$, respectively.
Notice that $\deg F^* = \sum_{i=0}^t [K_i:\K]$.

For an augmented $(K_0,\dots , K_t)$-form $F^*=(F,\theta_{0,F}, \dots ,\theta_{t,F})$
and for $A\in \text{GL}_2(\K)$, $\lambda \in \K^*$ we define
\begin{equation}\label{2.1}
\lambda F_A^*:=
(\lambda F_A,\langle A \rangle^{-1}\theta_{0,F},\dots ,\langle A \rangle^{-1}\theta_{t,F}).
\end{equation}
Clearly, $\lambda F_A^*$ is again an augmented $(K_0,\dots , K_t)$-form.
Notice that if $G^*=\lambda F_A^*$ then $F^*=\lambda^{-1} G_{A^{-1}}^*$;
further if $G^*=\lambda F_A^*$, $H^*=\mu G_B^*$ for some
$A,B \in \text{GL}_2(\K)$, $\lambda, \mu \in \K^*$ then
$H^*=\lambda \mu F_{AB}^*$.

Let $R$ be a subring of $\K$. Two augmented $(K_0, \dots , K_t)$-forms
$F^*, G^*$ are called $R$-\textit{equivalent}, notation $F^* \simm{R} G^*$, if
$G^*=F_U^*$ for some $U \in \text{GL}_2(R)$, and \textit{weakly}
$R$-\textit{equivalent},
notation $F^* \simmm{R} G^*$, if $G^*=\lambda F_U^*$ for some
$U \in \text{GL}_2(R)$ and $\lambda \in R^*$.

Let
$$
\text{M}_2^{\text{ns}}(R)=
\left\{ \begin{pmatrix} a&b\\c&d\end{pmatrix} : a,b,c,d \in R,\
\det \begin{pmatrix} a&b\\c&d\end{pmatrix}\ne 0 \right\}.
$$
Then for two augmented  $(K_0, \dots , K_t)$-forms $F^*, G^*$
we write $F^* \kacs{R} G^*$ if $G^* = F_A^*$ for some
$A \in \text{M}_2^{\text{ns}}(R)$.

In the Lemma below we have collected some simple facts.

\begin{lemma}\label{L2.1}
Let $r:=\sum_{i=0}^t[K_i:\K] \geq 3$ and let $R$ be a subring of $\K$.
%\begin{enumerate}
\vskip 0.1cm
%\item
{\bf (i)} Let $F^*$ be an augmented $(K_0, \dots , K_t)$-form,
$U \in \text{GL}_2(\K)$ and $\lambda \in \K^*$.
Then $\lambda F_U^*=F^*$ if and only if $U = \rho \bigl( \begin{smallmatrix} 1&0\\0&1\end{smallmatrix} \bigr)$
with $\rho \in \K^*$ and $\rho^r = \lambda^{-1}$.
\vskip 0.1cm

%\item
{\bf (ii)} Let $F^*, G^*$ be two augmented $(K_0, \dots , K_t)$-forms
and suppose that $G^*=\lambda_0 F_{U_0}^*$ for some $U_0 \in \text{GL}_2(\K)$,
$\lambda_0 \in \K^* $. Then for any other $U \in \text{GL}_2(\K)$,
$\lambda \in \K^*$ we have $G^*=\lambda F_U^*$ if and only if $U = \rho U_0$
with $\rho \in \K^*$ and $\rho^r = \lambda_0/\lambda$.
\vskip 0.1cm

%\item
{\bf (iii)} Let $F^*, G^*, H^*$ be augmented $(K_0, \dots , K_t)$-forms
such that
$F^* \kacs{R} G^*$, $G^* \kacs{R} H^*$. Then  $F^* \kacs{R} H^*$.
\vskip 0.1cm

%\item
{\bf (iv)}
Let $F^*, G^*$ be two augmented $(K_0, \dots , K_t)$-forms. Then\\
$F^* \kacs{R} G^*$, $\ G^* \kacs{R} F^* \iff F^* \simm{R} G^*$.
%\end{enumerate}
\end{lemma}

\begin{proof}
(i) Let $F^*=(F,\theta_{0,F}, \dots ,\theta_{t,F})$.
For $i= 0\kdots t$,
put $r_i:=[K_i:\K]$ and denote by
$\theta_{i,F}^{(1)}\kdots \theta_{i,F}^{(r_i)}$
the conjugates of $\theta_{i,F}$ over $\K$ (if
$\theta_{i,F}=\infty$, then $K_i=\K$, $r_i=1$ and
$\theta_{i,F}^{(1)}=\infty$). By assumption,
$\langle U \rangle^{-1}\theta_{i,F}=\theta_{i,F}$ for $i=0, \dots , t$ and
therefore, $\langle U \rangle^{-1}\theta_{i,F}^{(j)}=\theta_{i,F}^{(j)}$
for $i=0, \dots , t$, $j=1, \dots , r_i$. Thus, $\langle U \rangle$ has
$\sum_{i=0}^t[K_i:\K]=r \geq 3$ fixpoints. It follows that $\langle U \rangle$
is the identity on $\PP^1$, hence $U = \rho \bigl( \begin{smallmatrix} 1&0\\0&1\end{smallmatrix} \bigr)$
with $\rho \in \K^*$. Now since $\lambda F_U=F$, we have
$F(X,Y)=\lambda F(\rho X, \rho Y)=\lambda \rho^r F(X,Y)$, hence
$\rho^r=\lambda^{-1}$. Conversely, if
$U = \rho \bigl( \begin{smallmatrix} 1&0\\0&1\end{smallmatrix} \bigr)$
with $\rho^r=\lambda^{-1}$,
then clearly, $\lambda F_U^*=F^*$.

(ii) Let $G^*=\lambda F_U^*$.
Then $(\lambda_0 \lambda^{-1}) F_{U_0U^{-1}}^* = F^*$. Apply (i).

(iii) Obvious.

(iv) $\Leftarrow$ is clear. Assume  $F^* \kacs{R} G^*,\ G^* \kacs{R} F^*$.
Then there are $A,B \in \text{M}_2^{\text{ns}}(R)$ such that
$G^* = F_A^*$, $F^* = G_B^*$.
Thus $F^* = F_{AB}^*$. Hence by (i),
$AB=\rho \bigl( \begin{smallmatrix} 1&0\\0&1\end{smallmatrix} \bigr)$ with $\rho^r=1$.
Now $\rho \in R$ and $A^{-1}=\rho^{-1}B =\rho^{r-1}B \in \text{M}_2^{\text{ns}}(R)$.
So $A \in \text{GL}_2(R)$ and $F^* \simm{R} G^*$.
\end{proof}
\vskip 0.3cm

Let again $S$ be a finite subset of $M_{\K}$ containing
all infinite places. For $v \not\in S$
(i.e. $v \in M_{\K} \setminus S$) define the local ring
$\OO_v=\{x \in \K : |x|_v \leq 1 \}$. We need a few probably well-known
local-to-global results, relating (weak) $\OO_v$-equivalence
of augmented forms for $ v \not\in S$ to  $\OO_S$-equivalence.
We have inserted the proofs for lack of a good reference.

\begin{lemma}\label{L2.2}
Let $F^*, G^*$ be two  augmented $(K_0, \dots , K_t)$-forms such that
$F^*, G^*$ are $\OO_v$-equivalent for every $v \not\in S$. Then
$F^*, G^*$ are $\OO_S$-equivalent.
\end{lemma}

\begin{proof}
By assumption, for every $v \not\in S$ there is $U_v \in \text{GL}_2(\OO_v)$
such that $G^*=F_{U_v}^*$. By Lemma \ref{L2.1}, (ii) for $v \not\in S$
we have $U_v=\rho_v U_0$ where $U_0$ is one of the matrices $U_v$ $(v \not\in S)$,
and $\rho_v \in \K^*$, $\rho_v^r=1$. Then clearly,  $G^*=F_{U_0}^*$ and
$U_0 \in \text{GL}_2(\OO_v)$ for $v \not\in S$, so
$U_0 \in \text{GL}_2(\OO_S)$. Lemma \ref{L2.2} follows.
\end{proof}
\vskip 0.3cm

The following result is more involved.

\begin{lemma}\label{L2.3}
Let $\CC$ be a collection of augmented $(K_0, \dots , K_t)$-forms
such that for every pair $F^*, G^* \in \CC$ we have that $F^*, G^*$ are
weakly $\OO_v$-equivalent for every $v \not\in S$. Let $s:=\# S$.
Then $\CC$ is contained in the union of at most $r^s$ $\OO_S$-equivalence
classes if $r$ is odd, and in the union of at most $r^sh_2(\OO_S)$
$\OO_S$-equivalence classes if $r$ is even.
\end{lemma}

Before proving Lemma \ref{L2.3} we make some preparations.

If $R$ is a domain with quotient field $K$, then by a fractional
$R$-ideal, we mean
a subset $\A \ne \{ 0 \}$ of $K$
such that $\lambda \A$ is an ideal of $R$ for
some $\lambda \in K^*$.
For $v\not\in S$, denote by $\fp_v$ the prime ideal of $\OO_S$
corresponding to $v$, i.e., $\fp_v=\{ x\in \OO_S:\, |x|_v<1\}$,
and by
$\ordv$ the discrete valuation corresponding to $v$.
Thus, $[x]=\prod_{v\not\in S}\fp_v^{\ordv (x)}$ for $x\in \K^*$.

Let $F^*,G^*\in\CC$. Thus, for every $v\not\in S$ there are
$U_v\in\text{GL}_2(\OO_v)$, $\lambda_v\in\OO_v^*$ such that
$G^*=\lambda_v F^*_{U_v}$.
Choose any $U\in\text{GL}_2(\K )$, $\lambda\in \K^*$
such that $G^*=\lambda F^*_U$. Then by (ii) of Lemma \ref{L2.1}, for each
$v\not\in S$ there is a $\rho_v\in \K^*$ such that
\begin{equation}\label{2.2}
U_v=\rho_v U\, ,\quad\lambda_v=\rho_v^{-r}\lambda\, .
\end{equation}
Define the $\OO_S$-fractional ideal
\begin{equation}\label{2.3}
\fa (F^*,G^*) :=\prod_{v\not\in S} \fp_v^{\ordv (\rho_v)}\, .
\end{equation}
This is well-defined, since for all but finitely many $v\not\in S$
we have $\lambda\in\OO_v^*$, whence $\rho_v\in\OO_v^*$, whence
$\ordv (\rho_v)=0$. Let $\fA (F^* ,G^*)$ denote the ideal class
of $\fa (F^* ,G^*)$, that is, $\{ \mu\cdot\fa (F^*,G^*):\, \mu\in\K^*\}$.

The fractional ideal $\fa (F^* ,G^*)$ depends on the particular choice
of $U_v,\, \lambda_v$ $(v\not\in S)$, $U,\lambda$, but its ideal class
$\fA (F^*,G^*)$ does not.
Indeed, for $v\not\in S$, choose
$U_v'\in\text{GL}_2(\OO_v)$, $\lambda_v'\in\OO_v^*$
such that $G^*=\lambda_v'F^*_{U_v'}$ and then choose $U'\in\text{GL}_2(\K )$
and $\lambda'\in \K^*$ such that $G^*=\lambda ' F^*_{U'}$.
By (ii) of Lemma \ref{L2.1} there are $\rho_v'\in \K^*$ such that
$U_v'=\rho_v'U'$, $\lambda_v'={\rho_v'}^{-r}\lambda '$ for $v\not\in S$.
This gives rise to a fractional ideal
$\fa '(F^*,G^*)=\prod_{v\not\in S} \fp_v^{\ordv (\rho_v')}$.
Again by (ii) of Lemma \ref{L2.1}, there is $\mu\in \K^*$ such that
$U'=\mu U$ and $\lambda' =\mu^{-r}\lambda$. This implies for $v\not\in S$
that $U_v'=\rho_v'\mu\rho_v^{-1}U_v$, hence $\rho_v'\mu\rho_v^{-1}\in\OO_v^*$,
and so $\ordv(\rho_v')=\ordv(\rho_v)-\ordv(\mu)$.
Therefore, $\fa '(F^*,G^*)=\mu^{-1}\fa (F^*,G^*)$.

\begin{lemma}\label{L2.4}
{\bf (i)}
Let $F^*,G^*\in\CC$. Then $\fA (F^*,G^*)^{\gcd (r,2)}$ is the principal
ideal class.
\vskip 0.1cm

{\bf (ii)}
Let $F^*,G^*\in\CC$ and suppose that $\fA (F^*,G^*)$ is the principal
ideal class. Then $F^*$, $G^*$ are weakly $\OO_S$-equivalent.
\vskip 0.1cm

{\bf (iii)}
Let $F^*,G^*,H^*\in\CC$. Then $\fA (F^*,H^*)=\fA (F^*,G^*)\cdot\fA (G^*,H^*)$.
\end{lemma}

\begin{proof} (i) According to \eqref{2.2}
we have for $v\not\in S$, that
\begin{eqnarray*}
&&\ordv(\rho_v^2)=\ordv(\det U_v(\det U)^{-1})=\ordv ((\det U)^{-1}\, ,\\
&&\ordv(\rho_v^r)=\ordv(\lambda\lambda_v^{-1})=\ordv (\lambda )\, ,
\end{eqnarray*}
and so according to \eqref{2.3}, $\fa (F^*,G^*)^2=[\det U ]^{-1}$ and
$\fa (F^*,G^*)^r= [\lambda ]$, 
where $[a]$ denotes the $\OO_S$-fractional ideal
generated by $a$. This implies (i).

(ii) Let $\fa (F^* ,G^*)$ be given by \eqref{2.2}, \eqref{2.3}.
Then by our
assumption, $\fa (F^* ,G^*)=[\rho ]$ with $\rho\in\K^*$.
This implies $\rho\rho_v^{-1}\in\OO_v^*$ for $v\not\in S$. Put
$V:=\rho U$, $\mu :=\rho^{-r}U$. Then $G^*=\mu F^*_V$.
Further, by \eqref{2.2}, we have for $v\not\in S$,
that $U_v=\rho_v\rho^{-1}V$,
$\lambda_v=(\rho_v\rho^{-1})^{-r}\mu$, which implies
$V\in\text{GL}_2(\OO_v)$ and $\mu\in \OO_v^*$. Hence $V\in\text{GL}_2(\OO_S)$
and $\mu\in\OO_S^*$. Our assertion (ii) follows.

(iii) Straightforward computation.
\end{proof}
\vskip 0.3cm

\begin{proof}[Proof of Lemma \ref{L2.3}.] Fix $F^*\in\CC$. We subdivide
$\CC$ into classes such that two augmented forms $G_1^*,G_2^*\in\CC$ 
are in the
same class if and only if their corresponding ideal classes
$\fA (F^*,G_1^*)$, $\fA (F^*,G_2^*)$ coincide.
Let $F_1^*\kdots F_h^*$ be a full system of representatives
for these classes. Notice that by (i) of Lemma \ref{L2.4}, we have
$h\leq 1$ if $r$ is odd, and $h\leq h_2(\OO_S)$ if $r$ is even.

Fix $i\in\{ 1\kdots h\}$ and take any $G^*$ from the class
represented by $F_i^*$. According to (iii) of Lemma \ref{L2.4},
we have that $\fA (F_i^*,G^*)$ is the principal
ideal class. So by (ii) of Lemma \ref{L2.4}, there are $U\in\text{GL}_2(\OO_S)$
and $\ve\in\OO_S^*$ such that $G^*=\varepsilon (F_i^*)_U$.
The group $\OO_S^*$ is the direct product of $s=\# S$
cyclic groups, with generators $\varepsilon_1\kdots\varepsilon_s$, say.
So we may write
$\varepsilon =\varepsilon_1^{w_1}\cdots\varepsilon_s^{w_s}\eta^r$, with
$w_1\kdots w_r\in\{ 0\kdots r-1\}$ and $\eta\in\OO_S^*$. Consequently,
$G^*=\varepsilon_1^{w_1}\cdots\varepsilon_s^{w_s}(F_i^*)_{\eta U}$.

It follows that $\CC$ falls apart in at most $r^sh$ $\OO_S$-equivalence
classes, each represented by $\varepsilon_1^{w_1}\cdots\varepsilon_s^{w_s}F_i^*$
for certain $w_1\kdots w_s\in\{ 0\kdots r-1\}$, $i\in\{ 1\kdots h\}$.
Lemma \ref{L2.3} follows.
\end{proof}
\vskip 0.6cm

\section{From $\K$-equivalence classes to $\OO_S$-equivalence classes.}\label{S3}

\setcounter{equation}{0}

We keep the notation introduced in \S \S \ref{S1}-\ref{S2}.
Let $K_0, \dots , K_t$  be a sequence of finite extensions of $\K$.
Let $\CC$ be a set of augmented $(K_0, \dots , K_t)$-forms
which are all $\K$-equivalent to one another, and such that every
$F^* =(F,\theta_{0,F}\kdots \theta_{t,F})\in\CC$ satisfies $F \in \OO_S[X,Y]$
and \eqref{1.12}.
We will show that $\CC$ is
contained in finitely many $\OO_S$-equivalence classes
and estimate from above the number of these classes. We first
localize at a place $v \not\in S$, and estimate from above the
number of $\OO_v$-equivalence classes containing $\CC$.
Then we use Lemma \ref{L2.2}.
\vskip 0.15cm

Let $v \in M_{\K}$ be a finite place.
Denote by $\OO_v$ the local ring of $v$ and by $\fp_v$ the maximal
ideal of $\OO_v$, i.e.,
\[
\OO_v=\{ x\in\K :\, |x|_v\leq 1\},\quad
\fp_v=\{ x\in\K :\, |x|_v<1\}.
\]
Put $Nv :=\# (\OO_v/\fp_v )$.

Given a finite extension $L$ of $\K$,
we denote by $\OO_{L,v}$ the integral
closure of $\OO_v$ in $L$. The ring $\OO_{L,v}$ is a principal ideal
domain with finitely many prime ideals.
Further, it is a free $\OO_v$-module.
The $v$-discriminant ideal of $L/\K$ is given by
the ideal of $\OO_v$,
\begin{equation}\label{3.4}
\fd_{L/\K,v}=D_{L/\K}(\al_1, \dots ,\al_r)\cdot \OO_v\, ,
\end{equation}
where $\al_1, \dots , \al_r$ is any $\OO_v$-module basis of $\OO_{L,v}$.
This does not depend on the choice of $\al_1, \dots , \al_r$.

We will often denote the fractional $\OO_{L,v}$-ideal
generated by $a_1, \dots , a_m$ by $[a_1, \dots , a_m]$; from the
context it will always be clear in which field $L$ we are working.
Given a polynomial
$f \in L[X_1, \dots , X_m]$,
we denote by $[f]$ the fractional $\OO_{L,v}$-ideal
generated by the coefficients of $f$.
Then according to Gauss' Lemma,
\begin{equation}\label{3.3}
[fg]=[f][g] \ \ \ \text{for}\ \ f,g \in  L[X_1, \dots , X_m].
\end{equation}

Below we need some properties for resultants. The resultant
of two binary forms $F=a\prod_{i=1}^r (X-\alpha_iY)$,
$G=b\prod_{j=1}^s (X-\beta_jY)$ is given by
\begin{equation}\label{3.3a}
R(F,G)=a^sb^r\prod_{i=1}^r\prod_{j=1}^s (\alpha_i-\beta_j)\, .
\end{equation}
The resultant $R(F,G)$
is a polynomial
in the coefficients of
$F$ and $G$ with rational integral coefficients. 
It is homogeneous of degree $s$ in the coefficients of $F$
and homogeneous of degree $r$ in the coefficients of $G$.
For binary forms $F_0\kdots F_t$ we have
\begin{equation}\label{3.4a}
D(F)= \left( \prod_{i=0}^t D(F_i) \right)
\cdot \prod_{0\leq i<j \leq t} R(F_i, F_j)^2 .
\end{equation}

Now let $K_0, \dots , K_t$ be a sequence of finite extensions of $\K$.
Denote the normal closure over $\K$ of the
compositum $K_0 \dots K_t$ by $L$.
Put $r_i:=[K_i:\K]$ $(i=0, \dots ,t)$ and
$r:=r_0 + \dots + r_t$.  For $i=0, \dots ,t$ let
$\xi \mapsto \xi^{(i,j)}$ ($j=1, \dots ,r_i$) denote the
$\K$-isomorphic embeddings of $K_i$ into $L$.

We prove some properties for augmented $(K_0, \dots , K_t)$-forms.

\begin{lemma}\label{L3.1}
{\bf (i)} Let
$F^*=(F,\theta_{0,F}, \dots ,\theta_{t,F})$ be an augmented
$(K_0, \dots , K_t)$-form.
\vskip 0.1cm

{\bf (i)} Let $v\in M_{\K}$ be a finite place and suppose $F \in \OO_v[X,Y]$.
Then there is an ideal $\fc_v$ of $\OO_v$ such that
$$
D(F)\cdot \OO_v=\fc_v^2 \fd_{K_0/\K,v} \dots \fd_{K_t/\K,v}\, .
$$

{\bf (ii)} Suppose that $F \in \OO_S[X,Y]$. Then there is  an ideal $\fc$ of $\OO_S$
such that
$$
D(F)\cdot \OO_S =\fc^2 \fd_{K_0/\K,S} \dots \fd_{K_t/\K,S}\, .
$$
\end{lemma}

\begin{proof}
(ii) follows by applying (i) for every $v \not\in S$. We prove (i).
Since $\OO_v$ is a principal ideal domain we may write $F=F_0 F_1 \dots F_t$, where
$F_i^*=(F_i, \theta_{i,F})$ is an augmented $K_i$-form and $F_i \in \OO_v[X,Y]$
for $i=0\kdots t$.
In view of \eqref{3.4a} and since $R(F_i,F_j)\in\OO_v$
for all $i,j$,
it suffices to show that $D(F_i)\cdot\OO_v =\fc_{v,i}^2 \fd_{K_i/\K,v}$
for some ideal $\fc_{v,i}$ of $\OO_v$.

Write
$F_i(X,Y)=a_0X^{r_i}+a_1X^{r_i-1}Y+\dots + a_{r_i}Y^{r_i}$, and put
$\omega_1=1$, $\omega_2=a_0\theta_{i,F}$, $\omega_3=a_0\theta_{i,F}^2+a_1\theta_{i,F}$,
$\dots$, $\omega_{r_i}=a_0\theta_{i,F}^{r_i-1}+a_1\theta_{i,F}^{r_i-2}+
\dots + a_{r_i-2}\theta_{i,F}$. Let  $\{ \al_1, \dots , \al_{r_i} \}$ be an
$\OO_v$-basis of $\OO_{K_i,v}$. Then since
$\omega_1, \dots ,\omega_{r_i} \in \OO_{K_i,v}$ we have
$\omega_i=\sum_{j=1}^{r_i} \xi_{ij} \al_j$ with $\xi_{ij} \in \OO_v$.
Invoking \eqref{1.5a} we obtain
\begin{equation*}
\begin{split}
D(F_i)\cdot\OO_v&=D_{K_i/\K}(\omega_1, \dots ,\omega_{r_i})\cdot\OO_v\\
&=\det(\xi_{ij})^2 D_{K_i/\K}(\al_1,\dots ,\al_{r_i})\cdot\OO_v
= \det(\xi_{ij}) ^2\fd_{K_i/\K,v}.
\end{split}
\end{equation*}
Now Lemma \ref{L3.1} follows.
\end{proof}
\vskip 0.3cm

Let again
$F^*=(F,\theta_{0,F}, \dots ,\theta_{t,F})$ be an augmented
$(K_0, \dots , K_t)$-form. Henceforth we fix a finite place
$v\in M_{\K}$ and
assume that $F\in \OO_v[X,Y]$.
For $i=0, \dots , t$, choose  $\al_{i,F}$, $\be_{i,F}$ such that
\begin{equation}\label{3.5}
\begin{split}
&\al_{i,F},\be_{i,F}\in \OO_{K_i,v},\ \ \frac{\al_{i,F}}{\be_{i,F}}=\theta_{i,F}, \ \
[\al_{i,F},\be_{i,F}]=[1] \ \text{if}\ \ \theta_{i,F}\ne\infty,\\
&\al_{i,F}\in \OO_v^*, \be_{i,F}=0 \ \text{if}\ \ \theta_{i,F}=\infty;
\end{split}
\end{equation}
this is possible since $\OO_{K_i,v}$ is a principal ideal domain.
We may write
\begin{equation}\label{3.6}
  F=\ve_F \prod_{i=0}^t \prod_{j=1}^{r_i} (\be_{i,F}^{(i,j)}X - \al_{i,F}^{(i,j)}Y)\quad\mbox{with $\ve_F\in\OO_v$, $\ve_F\not=0$.}
\end{equation}
Indeed, a priori we know only that $\ve_F \in \K^*$.
But by Gauss' Lemma we have
\begin{equation}\label{3.7}
[F]=[\ve_F] \prod_{i=0}^t \prod_{j=1}^{r_i}
[\be_{i,F}^{(i,j)},\al_{i,F}^{(i,j)}]=[\ve_F],
\end{equation}
and thus $\ve_F \in \OO_v$ follows.

To pass from double to single indices we define a map
\begin{equation}\label{3.8}
\begin{split}
\vp:1,\dots ,r \ \to \  &
(0,1),\dots ,(0,r_0),\dots\\
&\ \ \dots , (1,1),\dots ,(1,r_1),\dots ,(t,1),\dots ,(t,r_t),
\end{split}
\end{equation}
meaning that $\vp$ maps $1,\dots ,r$ to $(0,1),\dots,(t,r_t)$,
respectively. We define the ideals of $\OO_{L,v}$:
\begin{equation}\label{3.9}
\fd_{kl}(F^*)=
[\al_{i_1,F}^{(i_1,j_1)}\be_{i_2,F}^{(i_2,j_2)}-
\al_{i_2,F}^{(i_2,j_2)}\be_{i_1,F}^{(i_1,j_1)}]
\end{equation}
for $k,l=1, \dots ,r$, $k<l$, where $\vp(k)=(i_1, j_1)$, $\vp(l)=(i_2,j_2)$.
Notice that the ideals $\fd_{kl}(F^*)$ are independent of the choice of
$\al_{i,F}, \be_{i,F}$ in \eqref{3.5}.
By \eqref{3.6}, \eqref{1.1},
we have
\begin{equation}\label{3.10}
\prod_{1\leq k<l \leq r } \fd_{kl}(F^*)^2 \supseteq [D(F)].
\end{equation}
Further, if  $G^*$ is an augmented $(K_0, \dots , K_t)$-form which is
$\OO_v$-equivalent to $F^*$ then
\begin{equation}\label{3.11}
\fd_{kl}(F^*)=\fd_{kl}(G^*) \ \ \text{for} \ \ 1 \leq k <l \leq r.
\end{equation}
The latter can be seen easily by taking
$U \in \text{GL}_2(\OO_v)$ such that $G^*=F_U^*$ and putting
$\bigl( \begin{smallmatrix} \al_{i,G}\\ \be_{i,G} \end{smallmatrix} \bigr):= U^{-1}
\bigl( \begin{smallmatrix} \al_{i,F}\\ \be_{i,F} \end{smallmatrix} \bigr)$,
$ \theta_{i,G} := \langle U \rangle^{-1}\theta_{i,F}$ for $i =0, \dots ,t$.
Then \eqref{3.5}, \eqref{3.6}, \eqref{3.9}
hold with everywhere $G,G^*$ in place
of $F,F^*$ and we obtain 
$\fd_{kl}(G^*)=(\det U^{-1})\cdot \fd_{kl}(F^*)=\fd_{kl}(F^*)$
since $\det U^{-1}\in \OO_v^*$.

\begin{lemma}\label{L3.2}
There are ideals $\fd_{kl}$ of $\OO_{L,v}$ independent of $F^*$ such that
\begin{equation}\label{3.12}
\fd_{kl}(F^*) \subseteq \fd_{kl} \ \ \text{for} \ \ 1 \leq k<l \leq r,
\end{equation}
\begin{equation}\label{3.13}
\prod_{1 \leq k<l \leq r} \fd_{kl}^2 \subseteq  \fd_{K_0/\K,v} \dots \fd_{K_t/\K,v} .
\end{equation}
\end{lemma}

\begin{proof}
Take $i \in \{ 0, \dots ,t\}$ and choose an $\OO_v$-basis
$\{ \al_{i,1}, \dots , \al_{i,r_i} \}$ of $\OO_{K_i,v}$.
Then there is a polynomial $I_{K_i/\K} \in \OO_v[X_1, \dots , X_{r_i}]$
(the index form of $K_i/\K$ with respect to
$ \al_{i,1}, \dots , \al_{i,r_i}$) such that
\begin{equation*}
\begin{split}
\prod_{1 \leq j_1<j_2 \leq r_i}
&\left( \sum_{m=1}^{r_i} \al_{i,m}^{(i,j_1)}X_m -
\sum_{m=1}^{r_i} \al_{i,m}^{(i,j_2)}X_m \right)^2 \\
& \qquad = D_{K_i/\K}(\al_{i,1}, \dots , \al_{i,r_i})
I_{K_i/\K}^2(X_1, \dots , X_{r_i}).
\end{split}
\end{equation*}
Define the ideal of $\OO_{L,v}$:
\begin{equation}\label{3.14}
\fb_{i,j_1,j_2}:=
\left[ \al_{i,1}^{(i,j_1)}- \al_{i,1}^{(i,j_2)}, \dots ,
\al_{i,r_i}^{(i,j_1)}- \al_{i,r_i}^{(i,j_2)} \right].
\end{equation}
Then by Gauss' Lemma
\begin{equation}\label{3.15}
\prod_{1 \leq j_1<j_2 \leq r_i} \fb_{i,j_1,j_2}^2 \subseteq
[D_{K_i/\K}(\al_{i,1}, \dots , \al_{i,r_i}) ] = \fd_{K_i/\K,v}.
\end{equation}
Moreover $\xi^{(i,j_1)} - \xi^{(i,j_2)} \in \fb_{i,j_1,j_2}$ for any
$\xi \in \OO_{K_i,v}$. Hence
for the numbers $\al_{i,F}, \be_{i,F}$ chosen in \eqref{3.9} we have
\begin{equation}\label{3.16}
\al_{i,F}^{(i,j_1)} \be_{i,F}^{(i,j_2)}-\al_{i,F}^{(i,j_2)} \be_{i,F}^{(i,j_1)}
\in \fb_{i,j_1,j_2} \ \ \ (1 \leq j_1<j_2 \leq r_i).
\end{equation}
Let $\vp$ be the map from \eqref{3.8}. Define $\fd_{kl}$ by
\begin{equation}\label{3.17}
\left\{
  \begin{aligned}
    &\fd_{kl}= \fb_{i,j_1,j_2} \ \ \ \ \text{if} \ \ \vp(k)=(i,j_1),\ \vp(l)=(i,j_2)\\
    &\fd_{kl}= [1] \ \ \ \ \ \ \ \ \text{if} \ \ \vp(k)=(i_1,j_1),\ \vp(l)=(i_2,j_2)\ \ \text{with} \ \ i_1 \ne i_2.
\end{aligned}
\right.
\end{equation}
Then \eqref{3.12}, \eqref{3.13} follow at once from \eqref{3.16}, \eqref{3.17},
\eqref{3.10}.
\end{proof}
\vskip 0.3cm

Let $\fc_v=\fc_v(F^*)$ be the ideal from (i) of Lemma \ref{L3.1}. Define
$\rho_v(F^*) \in \Z$ by $\fc_v=\fp_v^{\rho_v(F^*)}$. Thus,
$[D(F)]=\fp_v^{2\rho_v(F^*)} \prod_{i=0}^t \fd_{K_i/\K,v}$.

\begin{lemma}\label{L3.3}
 Let $\rho$ be a non-negative integer. Then as the tuple
$F^*=\\(F,\theta_{0,F}, \dots ,\theta_{t,F})$ runs
through the collection of augmented $(K_0, \dots , K_t)$-forms with
\begin{eqnarray}\label{3.18}
 &F\in \OO_v[X,Y]\, &
 \\[0.2cm]
\label{3.19}
&\rho_v(F^*) \leq \rho\, ,&
\end{eqnarray}
the tuple $(\fd_{kl}(F^*) : 1\leq k < l \leq r)$ runs through a set of
cardinality at most
\begin{equation}\label{3.20}
\binom{2\rho + \frac{1}{2}r(r-1)}{\frac{1}{2}r(r-1)}
\end{equation}
depending only on $K_0, \dots , K_t$, $v$, $\rho$.
\end{lemma}

\begin{proof}
We define an action of the Galois group $\Gal(L/\K)$ on the set
of subscripts $\{ 1\kdots r\}$ as follows.
Denote by $A$ the set of all $r$-tuples
$(\ga_1\kdots\ga_r)$ with the property that there are
$\xi_0\in K_0$, $\xi_1\in K_1\kdots \xi_t\in K_t$
such that
\begin{equation*}
(\ga_1\kdots\ga_r)=(\xi_0^{(0,1)},\dots , \xi_0^{(0,r_0)}, \dots ,
    \xi_t^{(t,1)},\dots , \xi_t^{(t,r_t)})\, .
\end{equation*}
Then there is a homomorphism $\tau\mapsto \tau^*$
from $\Gal(L/\K)$ to the permutation group of $\{ 1\kdots r\}$,
such that
\begin{equation}\label{3.22}
\tau(\ga_k)=\ga_{\tau^*(k)} \ \ \text{for}\ \ (\ga_1, \dots ,\ga_r) \in A, \ \ k=1, \dots ,r.
\end{equation}
Notice that if $\vp(k)=(i,j)$, then $\vp(\tau^*(k))=(i,j')$ for some
$j'\in \{1, \dots ,r_i \}$  where $\vp$ is the map given by \eqref{3.8}.

For each $k,l \in \{1, \dots ,r \} $, with $k<l$, we define the subfield
$L_{kl}$ of $L$ by
\begin{equation}\label{3.23}
  \Gal(L/L_{kl})=\left\{ \tau\in \Gal(L/\K) : \tau^*(\{k,l\})=\{k,l\} \right\}
\end{equation}
(i.e. $\tau^*(k)=k, \tau^*(l)=l$, or $\tau^*(k)=l, \tau^*(l)=k$).
We partition the set of pairs
$\{(k,l):\, k,l\in \{1, \dots , r\}, k<l \}$
into orbits $C_1, \dots , C_n$
in such a way that $(k_1,l_1)$, $(k_2,l_2)$
belong to the same orbit if and only if
$\{k_2,l_2\}=\tau^*(\{k_1,l_1\})$ for some $\tau \in \Gal(L/\K)$.
For each $m = 1, \dots ,n$ we choose a representative $(k_m,l_m)$ of $C_m$.
Then if
$(k,l)$ runs through $C_m$, the field
$L_{kl}$ runs through all conjugates over $\K$ of $L_{k_ml_m}$, and so
\begin{equation}\label{3.24}
\# C_m=[L_{k_ml_m}:\K ] \ \ \text{for} \ \ m = 1, \dots ,n.
\end{equation}

Now let $F^*=(F,\theta_{0,F}, \dots ,\theta_{t,F})$ be an augmented
$(K_0, \dots , K_t)$-form satisfying \eqref{3.18}, \eqref{3.19}. Define the ideals
\begin{equation*}
\A_{kl}(F^*):=\fd_{kl}(F^*)^2\fd_{kl}^{-2}\quad (1\leq k<l\leq r).
\end{equation*}
By Lemma \ref{L3.2} we have $\A_{kl}(F^*) \subseteq \OO_{L,v}$,
and by \eqref{3.9}, \eqref{3.14}, \eqref{3.17}, the ideal $\A_{kl}(F^*)$
is generated by elements from the field $L_{kl}$.
It is clear that the ideals $\A_{kl}(F^*)$ determine 
$\fd_{kl}(F^*)$ $(1 \leq k<l \leq r)$ uniquely.

For brevity put
\[
L_m :=L_{k_ml_m},\quad
\A_m(F^*):=\A_{k_ml_m}(F^*) \cap L_m\ \
\mbox{($m=1, \dots , n$);}
\]
thus $\A_m(F^*)$ is an ideal of $\OO_{L_m,v}$.
The ideals $\A_1(F^*), \dots ,\A_n(F^*)$ determine $\fd_{kl}(F^*)$
$(1 \leq k < l \leq r)$ uniquely. Indeed, they determine the ideals
$\A_{k_ml_m}(F^*)$ ($m=1, \dots ,n$) of $\OO_{L,v}$ since the latter are
generated by elements from $L_m$; and then by taking conjugates over $\K$ one
obtains all ideals $\A_{kl}(F^*)$ $(1 \leq k < l \leq r)$,
which, as mentioned before,
determine $\fd_{kl}(F^*)$ $(1 \leq k < l \leq r)$.

For $m=1, \dots ,n$ let $\fP_{m1}, \dots ,\fP_{mg_m}$
be the prime ideals of $\OO_{L_m,v}$. Thus,
\begin{equation*}
  \A_m(F^*)=\fP_{m1}^{w_{m1}(F^*)} \dots \fP_{mg_m}^{w_{mg_m}(F^*)}
\end{equation*}
where $w_{m1}(F^*), \dots ,w_{mg_m}(F^*)$ are non-negative integers since
$\A_m(F^*)$ is an ideal of $\OO_{L_m,v}$. Now the tuple of integers
\[
\underline{w}({F^*}) :=(w_{m,k}(F^*): m =1 , \dots , n, \ k=1, \dots , g_m)
\]
determines uniquely the ideals $\A_m(F^*)$ $(m =1 , \dots , n)$, hence the ideals
$\fd_{kl}(F^*)$ $(1 \leq k < l \leq r)$.
Therefore it suffices to show that for $\underline{w}({F^*})$
there are at most $\binom{2\rho + \frac{1}{2}r(r-1)}{\frac{1}{2}r(r-1)}$ possibilities.

Now on the one hand we have
by \eqref{3.10}, \eqref{3.13}, (i) of Lemma \ref{L3.1}, and assumption
\eqref{3.19},
\begin{equation*}
\begin{split}
\prod_{1 \leq k<l \leq r} \A_{kl}(F^*) &\supseteq
D(F) (\fd_{K_0/\K,v} \dots \fd_{K_t/\K,v})^{-1}\cdot\OO_{L,v}=\fc_v^2\OO_{L,v}
\\
&\supseteq \fp_v^{2\rho}\cdot\OO_{L,v}\, ,
\end{split}
\end{equation*}
while on the other hand,
\begin{equation*}
\begin{split}
\prod_{1 \leq k<l \leq r} \A_{kl}(F^*)
&=\prod_{m=1}^n\prod_{(k,l)\in C_m} \A_{kl}(F^*) =
\prod_{m=1}^n N_{L_m/\K}(\A_m(F^*))\cdot\OO_{L,v}
\\
&=\prod_{m=1}^n\prod_{h=1}^{g_m} N_{L_m/\K}(\fP_{mh})^{w_{mh}(F^*)}
\cdot \OO_{L,v}
\\
&=\prod_{m=1}^n\prod_{h=1}^{g_m} \fp_v^{f_{mh}w_{mh}(F^*)}\cdot\OO_{L,v}
\\
&\subseteq \fp_v^{ \sum_{m=1}^n\sum_{h=1}^{g_m} w_{mh}(F^*) }\cdot \OO_{L,v} ,
\end{split}
\end{equation*}
where $f_{mh}$ is the residue class degree of $\fP_{mh}$ over $\fp_v$.
Therefore,
\begin{equation}\label{3.28}
  \sum_{m=1}^n\sum_{h=1}^{g_m} w_{mh}(F^*) \leq 2 \rho.
\end{equation}
Now $g_m \leq [L_m:\K] \leq \#C_m$ for $m = 1 , \dots , n$
in view of \eqref{3.24}. 
Hence the number of
summands on the left-hand side is at most
$$
\sum_{m=1}^n \#C_m= \# \{(k,l): 1\leq k < l \leq r \} =\frac{1}{2}r(r-1).
$$
By elementary combinatorics,
the number of tuples of non-negative integers
$\underline{w}({F^*})$ with \eqref{3.28} is at most
$$
\binom{2\rho + \frac{1}{2}r(r-1)}{\frac{1}{2}r(r-1)}.
$$
As observed above, this implies Lemma \ref{L3.3}.
\end{proof}
\vskip 0.3cm

Let $\CC$ be a $\K$-equivalence class of augmented $(K_0, \dots , K_t)$-forms.
Given an ideal $\fc_v$ of $\OO_v$ and a tuple of ideals
$\{ \fd_{kl} : 1\leq k<l \leq r \}$ of $\OO_{L,v}$, let
$\CC(\fc_v, \{ \fd_{kl} \})$ denote the collection of augmented $(K_0, \dots , K_t)$-forms
$F^*=(F,\theta_{0,F}, \dots ,\theta_{t,F})$ such that
\begin{align}
&F^*  \in \CC; \label{3.28a} \\
&F \in \OO_v[X,Y];  \label{3.29} \\
&[D(F)]=\fc_v^2 \cdot \fd_{K_0/\K,v} \dots \fd_{K_t/\K,v}; \label{3.30} \\
&\fd_{kl}(F^*)=\fd_{kl} \ \ \text{for} \ \ k,l \in \{1, \dots ,r \}, \ 1 \leq k<l\leq r.
\label{3.31}
\end{align}

\begin{lemma}\label{L3.4}
Suppose $r:=\sum_{i=0}^t [K_i:\K] \geq 3$. Let $\fc_v$ be an ideal of $\OO_v$
and $\{ \fd_{kl} : 1\leq k<l \leq r \}$ a collection of ideals from $\OO_{L,v}$
such that the set $\CC(\fc_v, \{ \fd_{kl} \})$ is not contained in a single
$\OO_v$-equivalence class. Then
\begin{equation}\label{3.31a}
\fc_v \subseteq \fp_v^{\frac{r(r-1)}{2}}, \ \
\fd_{kl} \subseteq \fp_v\OO_{L,v} \ \ \text{for}\ \
1 \leq k<l \leq r,
\end{equation}
and for every $F^* \in \CC(\fc_v, \{ \fd_{kl} \})$
there is an $H^*$ with
\begin{equation}\label{3.31b}
H^* \kacs{\OO_v} F^*\, ,\quad H^* \in
\CC(\fp_v^{-\frac{1}{2}r(r-1)} \fc_v,\{ \fp_v^{-1}\fd_{kl} \})\, .
\end{equation}
\end{lemma}

\begin{proof}
If $H^*=(H,\theta_{0,H}\kdots \theta_{t,H})$ is an augmented form
with $H\in\OO_v[X,Y]$, then $\fd_{kl}(H^*)$ ($1\leq k<l\leq r$) are
all ideals of $\OO_{L,v}$, and by (i) of Lemma \ref{L3.1},
there is an ideal $\fc_v'\subseteq \OO_v$ such that
$[D(H)]=\fc_v'^2 \cdot \fd_{K_0/\K,v} \dots \fd_{K_t/\K,v}$.
So if we have shown that
there exists an $H^*$ with \eqref{3.31b},
then \eqref{3.31a} follows automatically.

Let  $F^* \in \CC(\fc_v, \{ \fd_{kl} \})$. There is a
$G^* \in \CC(\fc_v, \{ \fd_{kl} \})$ which is not $\OO_v$-equivalent to $F^*$.
This means that there is a matrix $A\in \text{GL}_2(\K)$ with
$A\not\in\text{GL}_2(\OO_v)$ such that $G^*=F_A^*$.
Since $\OO_v$ is a principal ideal domain, there are matrices
$U_1,U_2 \in \text{GL}_2(\OO_v)$ such that
\begin{equation*}
A=U_1  \begin{pmatrix} \al&0\\0&\de \end{pmatrix} U_2
\end{equation*}
with
\begin{equation}\label{3.33}
\al,\de \in \K^*, \quad  \frac{\de}{\al} \in \OO_v,\quad
\begin{pmatrix} \al&0\\0&\de \end{pmatrix} \not\in \text{GL}_2(\OO_v).
\end{equation}

Put $\til{F}^*:= F^*_{U_1}$, $\til{G}^*:= G^*_{U_2^{-1}}$.
Then
\begin{equation}\label{3.34}
\til{G}^*=\til{F}_{\bigl( \begin{smallmatrix} \al&0\\0&\de \end{smallmatrix} \bigr)}^*\, .
\end{equation}
Further,
$\til{F}^* \simm{\OO_v} F^*$, $\til{G}^* \simm{\OO_v} G^*$,
so by \eqref{3.11}, \eqref{1.4},
\begin{equation}\label{3.36}
\til{F}^*, \til{G}^* \in \CC(\fc_v, \{ \fd_{kl} \}).
\end{equation}
Clearly, in view of (iv) of Lemma \ref{L2.1}, it follows that there
is an $H^*$ with \eqref{3.31b} once we have proved
that there is an $H^*$ with
\begin{equation}\label{3.35}
H^*\kacs{\OO_v} \til{F}^*,\quad
H^* \in \CC(\fp_v^{-\frac{1}{2}r(r-1)} \fc_v,\{ \fp_v^{-1}\fd_{kl} \})\, .
\end{equation}

By \eqref{3.36},\eqref{3.30}, \eqref{3.34}, \eqref{1.3}, 
we have
\begin{equation*}
[D(\til{F})]=\fc_v^2 \prod_{i=0}^t \fd_{K_i/\K,v}=[D(\til{G})]=[\al\de]^{r(r-1)}[D(\til{F})]\, ,
\end{equation*}
and together with \eqref{3.33} this implies
\begin{equation}\label{3.37}
\de \in \OO_v, \quad \de \not\in \OO_v^*,\quad \al\de \in \OO_v^*.
\end{equation}

Write $\til{F}^*=(\til{F},\theta_{0,\til{F}}, \dots ,\theta_{t,\til{F}})$.
Then by
\eqref{3.34} we have
\[
\til{G}^*=\Big(\til{F}_{\bigl( \begin{smallmatrix} \al&0\\0&\de \end{smallmatrix} \bigr)},
\frac{\de}{\al}\theta_{0,\til{F}}, \dots ,\frac{\de}{\al}\theta_{t,\til{F}}
\Big)\, .
\]
Similarly as in \eqref{3.5},
choose $\al_{i,\til{F}},\be_{i,\til{F}} \in \OO_{K_i,v}$
such that $\al_{i,\til{F}}/\be_{i,\til{F}}=\theta_{i,\til{F}}$ and
$[\al_{i,\til{F}},\be_{i,\til{F}}]=[1]$ if $\theta_{i,\til{F}} \ne \infty$, and
$\al_{i,\til{F}} \in \OO_v^*$, $\be_{i,\til{F}}=0$ if $\theta_{i,\til{F}}=\infty$.
Likewise, choose $\al_{i,\til{G}},\be_{i,\til{G}} \in \OO_{K_i,v}$
such that $\al_{i,\til{G}}/\be_{i,\til{G}}=\de \theta_{i,\til{F}}/\al$ and
$[\al_{i,\til{G}},\be_{i,\til{G}}]=[1]$ if $\theta_{i,\til{F}} \ne \infty$, and
$\al_{i,\til{G}} \in \OO_v^*$, $\be_{i,\til{G}}=0$ if $\theta_{i,\til{F}}=\infty$.
Then  for $i=0, \dots ,t$ there is a $\lambda_i \in K_i^*$ such that
\begin{equation}\label{3.38}
(\al_{i,\til{G}},\be_{i,\til{G}} ) = \lambda_i (\de \al_{i,\til{F}},\al \be_{i,\til{F}} )
\ \ \text{for} \ \ i=0, \dots ,t.
\end{equation}

Take two pairs $(i_1,j_1)$, $(i_2,j_2)$ from
$\{ (i,j) : i=0, \dots ,t,\, j=1, \dots ,r_i \}$. Let
$k,l\in \{ 1, \dots ,r \}$ be such that
$\vp(k)=(i_1,j_1), \vp(l)=(i_2,j_2)$, where $\vp$ is the map from \eqref{3.8}.
Then by \eqref{3.36}, \eqref{3.38}, \eqref{3.37} and again \eqref{3.36},
\begin{equation*}
  \begin{split}
    \fd_{kl}=&
[\al_{i_1,\til{G}}^{(i_1,j_1)}\be_{i_2,\til{G}}^{(i_2,j_2)}-
\al_{i_2,\til{G}}^{(i_2,j_2)}\be_{i_1,\til{G}}^{(i_1,j_1)}] \\
=& [\lambda_{i_1}^{(i_1,j_1)} \lambda_{i_2}^{(i_2,j_2)}\al\de
(\al_{i_1,\til{F}}^{(i_1,j_1)}\be_{i_2,\til{F}}^{(i_2,j_2)}-
\al_{i_2,\til{F}}^{(i_2,j_2)}\be_{i_1,\til{F}}^{(i_1,j_1)})] \\
=&
[\lambda_{i_1}^{(i_1,j_1)}] [\lambda_{i_2}^{(i_2,j_2)}] \fd_{kl}
  \end{split}
\end{equation*}
and so $[\lambda_{i_1}^{(i_1,j_1)}] [\lambda_{i_2}^{(i_2,j_2)}]=[1]$.
This holds for any two distinct pairs $(i_1,j_1)$, $(i_2,j_2)$ from
$\{ (i,j) : i=0, \dots ,t,\, j=1, \dots ,r_i \}$.
Taking any pair $(i,j)$ from this set and then any two other pairs
$(i_1,j_1)$, $(i_2,j_2)$ (which is possible since
by assumption $r_0+\dots+r_t =r \geq 3$),
we obtain
$$
[\lambda_{i}^{(i,j)}]^2=
\frac{[\lambda_{i}^{(i,j)}] [\lambda_{i_1}^{(i_1,j_1)}] [\lambda_{i}^{(i,j)}]
[\lambda_{i_2}^{(i_2,j_2)}]}{[\lambda_{i_1}^{(i_1,j_1)}][\lambda_{i_2}^{(i_2,j_2)}] }=[1],
$$
so $[\lambda_{i}^{(i,j)}]=[1]$ for $i = 0, \dots ,t$, $j=1, \dots , r_i$. Together with \eqref{3.38}, this implies 
\[
[\de \al_{i,\til{F}},\al \be_{i,\til{F}} ]=[1]
\quad\mbox{for $i = 0, \dots ,t$.}
\]

By \eqref{3.37} we have $\de\in\fp_v$, hence
$\de\al_{i,\til{F}}\in\fp_v\OO_{L,v}$
for $i = 0, \dots ,t$. This implies that $\de\al_{i,\til{F}}$ is divisible
by each prime ideal of $\OO_{L,v}$,
therefore $[\al \be_{i,\til{F}}]=[1]$ for $i = 0, \dots ,t$.
Since by \eqref{3.37}, $[\al]=[\de^{-1}] \supseteq \fp_v^{-1}$ we have
$\be_{i,\til{F}}\in\fp_v\OO_{L,v}$  for $i = 0, \dots ,t$. So
\begin{equation}\label{3.39}
\be_{i,\til{F}}^{(i,j)}\in \fp_v\OO_{L,v}\ \
\text{for $i = 0, \dots ,t$, $j=1, \dots , r_i$}.
\end{equation}

We now construct an $H^*$ with \eqref{3.35}.
Choose $\Pi$ with $\fp_v=[\Pi ]$ and take
$$
H^*=\til{F}_{\bigl( \begin{smallmatrix} \Pi^{-1}&0\\0&1 \end{smallmatrix} \bigr)}^*=
(\til{F}_{\bigl( \begin{smallmatrix} \Pi^{-1}&0\\0&1 \end{smallmatrix} \bigr)},
\Pi\theta_{0,\til{F}}, \dots , \Pi\theta_{t,\til{F}}).
$$
Clearly,
\begin{equation}\label{3.39a}
H^* \kacs{\OO_v} \til{F}^*\, .
\end{equation}

Similarly as in \eqref{3.6} we may write
$$
 \til{F}=\ve_{\til{F}} \prod_{i=0}^t \prod_{j=1}^{r_i}
(\be_{i,\til{F}}^{(i,j)}X - \al_{i,\til{F}}^{(i,j)}Y)
\ \ \text{with} \ \ \ve_{\til{F}} \in \OO_v.
$$
Now \eqref{3.39} implies that
$$
H:=\til{F}_{\bigl( \begin{smallmatrix} \Pi^{-1}&0\\0&1 \end{smallmatrix} \bigr)}=
\ve_{\til{F}} \prod_{i=0}^t \prod_{j=1}^{r_i}
(\Pi^{-1}\be_{i,\til{F}}^{(i,j)}X - \al_{i,\til{F}}^{(i,j)}Y) \in \OO_{L,v}[X,Y].
$$
Since also $H \in \K[X,Y]$, we have
\begin{equation}\label{3.40}
H \in \OO_v[X,Y].
\end{equation}
Moreover, by \eqref{1.3}, \eqref{3.36},
\begin{equation}\label{3.41}
[D(H)]=[\Pi^{-r(r-1)}D(\til{F}) ] =
(\fp_v^{-\frac{1}{2}r(r-1)} \fc_v)^2 \fd_{K_0/\K,v} \dots \fd_{K_t/\K,v}.
\end{equation}
Further, we have $\Pi\theta_{i,\til{F}} = \al_{i,\til{F}}/\Pi^{-1}  \be_{i,\til{F}}$
and $[\al_{i,\til{F}}, \Pi^{-1}  \be_{i,\til{F}}]=[1]$ for $i = 0, \dots ,t$.
The latter is true since $\al_{i,\til{F}}, \Pi^{-1}  \be_{i,\til{F}} \in \OO_{L,v}$
and $[\al_{i,\til{F}},  \be_{i,\til{F}}]=[1]$. So by definition \eqref{3.9} and
by \eqref{3.36} we have for $1 \leq k<l \leq r$,
\begin{equation}\label{3.42}
\begin{split}
\fd_{kl}(H^*)=&[ \al_{i_1,\til{F}}^{(i_1,j_1)}\Pi^{-1}\be_{i_2,\til{F}}^{(i_2,j_2)}-
\al_{i_2,\til{F}}^{(i_2,j_2)}\Pi^{-1}\be_{i_1,\til{F}}^{(i_1,j_1)} ]\\
=&[\Pi]^{-1}\fd_{kl}(\til{F}^*)=\fp_v^{-1}\fd_{kl},
\end{split}
\end{equation}
where $\vp(k)=(i_1,j_1), \vp(l)=(i_2,j_2)$.

Now by collecting \eqref{3.39a}, \eqref{3.40}, \eqref{3.41}, \eqref{3.42}
and the obvious fact that $H^*$ is $\K$-equivalent to $\til{F}^*$
we infer that indeed $H^*$ satisfies \eqref{3.35}.
This completes the proof of Lemma \ref{L3.4}.
\end{proof}
\vskip 0.3cm

\begin{lemma}\label{L3.5}
Suppose $r:=\sum_{i=0}^t [K_i:\K] \geq 3$. Let
$\fc_v, \{ \fd_{kl} : 1\leq k<l \leq r \}$ be as in Lemma \ref{L3.4}.
Suppose that $\CC(\fc_v, \{ \fd_{kl} \}) \ne \emptyset$. Then there is
an augmented $(K_0, \dots , K_t)$-form
$F_0^*=(F_0,\theta_{0,F_0}, \dots ,\theta_{t,F_0})$ such that
\begin{equation}\label{3.43}
F_0 \in \OO_v[X,Y]
\end{equation}
and
\begin{equation}\label{3.44}
  F_0^* \kacs{\OO_v} F^* \ \ \text{for every} \ \ F^*\in \CC(\fc_v, \{ \fd_{kl} \}).
\end{equation}
\end{lemma}

\begin{proof}
We claim that there is a non-negative integer $i$ such that
\begin{eqnarray}
\label{3.45}
&&\fp_v^{-\frac{1}{2}r(r-1)i}\fc_v \subseteq [1], \ \
\fp_v^{-i}\fd_{kl}  \subseteq [1]
\ \ (1 \leq k<l \leq r)\, ,
\\
\label{3.46}
&&\CC(\fp_v^{-\frac{1}{2}r(r-1)i} \fc_v,\{ \fp_v^{-i}\fd_{kl} \}) \ne \emptyset\, ,
\\
\label{3.47}
&&
\CC(\fp_v^{-\frac{1}{2}r(r-1)i} \fc_v,\{ \fp_v^{-i}\fd_{kl} \})\\
\nonumber
&&\qquad\qquad\text{is contained in a single $\OO_v$-equivalence class.}
\end{eqnarray}

Indeed, if there is no such integer $i$, then by inductively
applying Lemma \ref{3.4} it follows that there are
arbitrarily large integers $i$ with \eqref{3.45}, \eqref{3.46}.
But there cannot be
arbitrarily large $i$ with \eqref{3.45}.

Let $i_0$ be the smallest integer $i$ with
\eqref{3.45}, \eqref{3.46}, \eqref{3.47}.
Pick 
\[
F_0^*=(F_0,\theta_{0,F_0}, \dots ,\theta_{t,F_0})
\in \CC(\fp_v^{-\frac{1}{2}r(r-1)i_0} \fc_v,\{ \fp_v^{-i_0}\fd_{kl} \})\, .
\]
Then $F_0[X,Y]\in  \OO_v[X,Y]$. By Lemma \ref{L3.4},
for every $F^*\in \CC(\fc_v, \{ \fd_{kl} \})$ there is a sequence
$$
F_{i_0}^* \kacs{\OO_v} F_{i_0-1}^* \kacs{\OO_v} \dots \kacs{\OO_v}
F_1^* \kacs{\OO_v} F^*
$$
with
$F_i^* \in \CC(\fp_v^{-\frac{1}{2}r(r-1)i} \fc_v,\{ \fp_v^{-i}\fd_{kl} \})$
for $i = 1, \dots ,i_0$.
By \eqref{3.47} we have $F_0^* \simm{\OO_v} F_{i_0}^*$
and then by (iv) and (iii) of Lemma \ref{L2.1},
$F_0^* \kacs{\OO_v} F_{i_0}^*$, $F_0^*\kacs{\OO_v} F^*$.
This proves Lemma \ref{L3.5}.
\end{proof}
\vskip 0.3cm

\begin{lemma}\label{L3.6}
Suppose $r:=\sum_{i=0}^t [K_i:\K] \geq 3$. Let $\fc_v$ be an ideal of $\OO_v$.
Let $\rho_v$ be the non-negative integer given by
$\fc_v=\fp_v^{\rho_v}$.
Let $\CC$ be a
$\K$-equivalence class of  augmented $(K_0, \dots , K_t)$-forms.
Denote by $\CC(\fc_v)$ the collection of  augmented $(K_0, \dots , K_t)$-forms
$F^*=(F,\theta_{0,F}, \dots ,\theta_{t,F})$ in $\CC$ satisfying
\begin{align}
 &F \in \OO_v[X,Y], \label{3.48}\\
 &[D(F)]=\fc_v^2 \fd_{K_0/\K,v} \dots \fd_{K_t/\K,v}. \label{3.49}
\end{align}
Then $\CC(\fc_v)$ is the union of at most
\begin{equation}\label{3.50}
  \binom{2\rho_v + \frac{1}{2}r(r-1)}{\frac{1}{2}r(r-1)}
\left( \sum_{i=0}^{[2\rho_v/r(r-1)]} (Nv)^i \right)
\end{equation}
$\OO_v$-equivalence classes.
\end{lemma}

\begin{proof}
By Lemma \ref{L3.3}, we can express the set $\CC(\fc_v)$ as a union of at most
$\binom{2\rho_v + \frac{1}{2}r(r-1)}{\frac{1}{2}r(r-1)}$ sets
$\CC(\fc_v, \{ \fd_{kl} \})$ where $\fd_{kl}$ $(1\leq k < l \leq r)$ are ideals
of $\OO_{L,v}$. So it suffices to show that for given ideals $\fc_v$ of $\OO_v$
and $\fd_{kl}$ $(1\leq k < l \leq r)$ of $\OO_{L,v}$, the set
$\CC(\fc_v, \{ \fd_{kl} \})$ is the union of not more than
\begin{equation}\label{3.51}
\sum_{i=0}^{[2\rho_v/r(r-1)]} (Nv)^i
\end{equation}
$\OO_v$-equivalence classes.

According to Lemma \ref{L3.5}, there is a fixed augmented $(K_0, \dots , K_t)$-form
$F_0^*=(F_0,\theta_{0,F}, \dots ,\theta_{t,F})$ with $F_0 \in \OO_v[X,Y]$
such that $F_0^* \kacs{\OO_v} F^*$ for every $F^* \in \CC(\fc_v, \{ \fd_{kl} \})$.
That is, for every $F^* \in \CC(\fc_v, \{ \fd_{kl} \})$ there is a matrix
$A \in \text{M}_2^{\text{ns}}(\OO_v)$ such that
$F^* = (F_0^*)_A$.
By Lemma \ref{L3.1}, there is an ideal $\fc_{v0}$ of $\OO_v$
such that $[D(F_0)]=\fc_{v0}^2 \fd_{K_0/\K,v} \dots \fd_{K_t/\K,v}$. Let
$\rho_{v0} \in \Z_{\geq 0}$ be defined by $\fc_{v0}=\fp_v^{\rho_{v0}}$.
Then by \eqref{3.49}, \eqref{1.3},
\begin{equation*}
\begin{split}
[D(F)]&=\fp_v^{2\rho_v} \fd_{K_0/\K,v} \dots \fd_{K_t/\K,v} \\
&=[\det A]^{r(r-1)}[D(F_0)]=
[\det A]^{r(r-1)}\fp_v^{2\rho_{v_0}} \fd_{K_0/\K,v} \dots \fd_{K_t/\K,v}.
\end{split}
\end{equation*}
Hence
\begin{equation}\label{3.52}
[\det A]=\fp_v^u \ \ \text{with} \ \ u=\frac{2(\rho_v-\rho_{v0})}{r(r-1)}.
\end{equation}
Choose $\Pi$ with $\fp_v =[\Pi ]$.
The ideals of $\OO_v$ are of the shape $\fp_v^m$ ($m\geq 0$) and
$\# \OO_v/\fp_v^m$ has cardinality $(Nv)^m$. From these facts it can be deduced
that every matrix $A \in \text{M}_2^{\text{ns}}(\OO_v)$ with \eqref{3.52}
can be expressed as
$$
A=A_{ij}U \ \ \text{with} \ \ U \in \text{GL}_2(\OO_v), \ \ \
A_{ij} = \begin{pmatrix}\Pi^{u-i} &0\\ \be_{ij} &\Pi^i  \end{pmatrix}
$$
where $i \in \{0,1, \dots ,u \}$ and where $\be_{i1},\dots ,\be_{i,(Nv)^i}$ is a
full system of representatives for the residue classes of $\OO_v$ modulo $\fp_v^i$.

Now if $F^* \in \CC(\fc_v, \{ \fd_{kl} \})$ then
$F^* = (F_0^*)_A$ for some $A \in \text{M}_2^{\text{ns}}(\OO_v)$ with
\eqref{3.52}, hence
$F^* = (F_0^*)_{A_{ij}U} \simm{\OO_v}(F_0^*)_{A_{ij}}$ for some
$i \in \{0, \dots ,u \}$, $j \in \{1, \dots ,(Nv)^i \}$. This implies that
$\CC(\fc_v, \{ \fd_{kl} \})$ is contained in the union of
$$
\sum_{i=0}^u (Nv)^i = \sum_{i=0}^{2(\rho_v-\rho_{v0})/r(r-1)} (Nv)^i
\leq \sum_{i=0}^{[2\rho_v/r(r-1)]} (Nv)^i
$$
$\OO_v$-equivalence classes. This proves Lemma \ref{L3.6}.
\end{proof}
\vskip 0.3cm
We now arrive at the main result of this section. We have formulated it
both for augmented forms and for ordinary binary forms.

\begin{proposition}\label{P3.7}
Let $\fc$ be an ideal of $\OO_S$. 
Let $r:=\sum_{i=0}^t [K_i:\K] \geq 3$.
\vskip 0.1cm

{\bf (i)} Let $\CC (\fc )$ be a $\K$-equivalence class of augmented 
$(K_0, \dots , K_t)$-forms such that any two elements of $\CC (\fc )$
are $\K$-equivalent and such that every
$F^*=(F,\theta_{0,F}, \dots ,\theta_{t,F})\in\CC (\fc )$ satisfies
\begin{align}
 &F \in \OO_S[X,Y]\, , \label{3.53} \\
 &D(F)\cdot\OO_S =  \fc^2 \fd_{K_0/\K,S} \dots \fd_{K_t/\K,S}\, . \label{3.54}
\end{align}
Then $\CC (\fc)$ is contained in
the union of at most
\begin{equation}\label{3.100}
\tau_{\frac{1}{2}r(r-1)}(\fc^2)
\left( \sum_{\fd^{\frac{1}{2}r(r-1)} \mid \fc} N_S(\fd) \right)
\end{equation}
$\OO_S$-equivalence classes.
\vskip 0.1cm

{\bf (ii)} Let $\CCC (\fc )$ be a subset of $\F (\OO_S ,K_0\kdots K_t)$
such that any two binary forms in $\CCC (\fc )$ are $\K$-equivalent
and such that every $F\in\CCC (\fc )$ satisfies \eqref{3.54}.
Then $\CCC (\fc )$ is contained in the union of finitely many 
$\OO_S$-equivalence classes, the number of which is bounded above
by \eqref{3.100}.
\end{proposition}

\begin{proof}
(i) For $v \not\in S$, let $\fp_v$ be the prime ideal of $\OO_S$ corresponding to $v$,
i.e., $\fp_v = \{ x \in \OO_S : |x|_v < 1 \}$. Then
$\fc = \prod_{v\not\in S} \fp_v^{\rho_v}$ with $\rho_v \in \Z_{\geq 0}$.
According to Lemma \ref{L3.6}, for each $v \not\in S$ the collection
$\CC  (\fc)$ is contained in the union of at most
\begin{equation*}
\begin{split}
A_v :=& \binom{2\rho_v + \frac{1}{2}r(r-1)}{\frac{1}{2}r(r-1)}
\sum_{i=0}^{[2\rho_v/r(r-1)]} (Nv)^i\\
=&
\binom{2\rho_v + \frac{1}{2}r(r-1)}{\frac{1}{2}r(r-1)}
\sum_{i=0}^{[2\rho_v/r(r-1)]} (N_S\fp_v)^i
\end{split}
\end{equation*}
$\OO_v$-equivalence classes. Lemma \ref{L2.2} implies that if
$\mathcal{A}_v$ is an $\OO_v$-equivalence class of
augmented $(K_0, \dots , K_t)$-forms for $v \not\in S$, then
$\cap_{v \not\in S} \mathcal{A}_v$ is an $\OO_S$-equivalence class.
This implies that $\CC  (\fc)$ is contained in the union of at most
$$
\prod_{v \not\in S} A_v = \tau_{\frac{1}{2}r(r-1)}(\fc^2)
\left( \sum_{\fd^{\frac{1}{2}r(r-1)} \mid \fc} N_S(\fd) \right)
$$
$\OO_S$-equivalence classes. This proves (i).

(ii) Fix $F_0\in\CCC (\fc )$. Extend $F_0$ to an augmented
$(K_0\kdots K_t)$-form $F_0^*=(F_0,\theta_{0,F_0}\kdots\theta_{t,F_t})$.
For every $F\in \CCC (\fc )$, choose $A\in\text{GL}_2(K)$ such that
$F=(F_0)_A$ and define $F^*:= (F_0^*)_A$. 
Clearly, the augmented forms constructed in this manner are $\K$-equivalent
to one another.
Now by applying (i) to the
collection $\CC (\fc):= \{ F^* :\, F\in\CCC (\fc )\}$, our assertion (ii)
follows at once.
\end{proof} 
\vskip 0.5cm

\section{Orders}\label{S4}

\setcounter{equation}{0}

Below, $\K$ is a number field, and $K$ is a finite extension of $\K$
of degree $r \geq 3$.
We denote by $\xi \mapsto \xi^{(i)}$ $(i=1, \dots , r)$
the $\K$-isomorphic embeddings of $K$ into some normal
closure $L$ of $K$ over $\K$.
As before, $S$ is a finite subset of $M_{\K}$ containing
all infinite places. Denote by $\OO_{L,S}$ the
integral closure of $\OO_S$ in $L$. Given $a_1, \dots , a_m$, we denote by
$[a_1, \dots , a_m]$ the fractional $\OO_{L,S}$-ideal
generated by $a_1, \dots , a_m$. For $f \in L[X_1, \dots , X_m]$
denote by $[f]$ the fractional $\OO_{L,S}$-ideal
generated by the coefficients of $f$. Given fractional
$\OO_{L,S}$-ideals $\fa$, $\fb$ we write $\frac{\fa}{\fb}$ for
$\fa \fb^{-1}$ where $\fb^{-1}$ is the inverse fractional $\OO_{L,S}$-ideal
of $\fb$. For a finitely generated $\OO_S$-module $\M \subset K$ with
$\M \ne (0)$ define
\begin{equation}\label{4.1}
\fd_{ij}(\M):=[\xi^{(i)}-\xi^{(j)} : \xi \in \M] \ \ \ (1 \leq i,j \leq r, \ i \ne j)
\end{equation}
to be the fractional $\OO_{L,S}$-ideal generated by all elements $\xi^{(i)}-\xi^{(j)}$
$(1 \leq i,j \leq r, \ i \ne j)$ with $\xi\in\M$ and
\begin{equation}\label{4.2}
\fD(\M):=[D_{K/\K}(\omega_1, \dots ,\omega_r) : \omega_1, \dots ,\omega_r \in \M ]
\end{equation}
to be the fractional $\OO_{L,S}$-ideal generated by all discriminants of all
$r$-tuples $\omega_1, \dots ,\omega_r \in \M$.

Let $F^*=(F,\theta_F)$ be an augmented $K$-form. Suppose that $F\in R[X,Y]$
where $R$ is some subring of $\K$.
Then the invariant order $\OO_{F^*,R}$
of $F^*$ is defined to be the $R$-submodule of $K$
with basis $\omega_1\kdots\omega_r$ given
by \eqref{1.5}.
By Simon \cite{Simon},
$\OO_{F^*,R}$ is indeed an $R$-order with quotient field $K$,
\begin{equation}\label{4.5}
F^* \simm{R} G^* \Rightarrow \OO_{F^*,R}=\OO_{G^*,R}
\end{equation}
for any two augmented $K$-forms $F^*,G^*$
(which is slightly stronger than \eqref{1.5b}),
and $D_{K/\K}(\omega_1, \dots ,\omega_r)=D(F^*)$. If $R=\OO_S$ we write
$\OO_{F^*,S}$ for $\OO_{F^*,R}$ and if $R=\OO_v$ (local ring) we write
$\OO_{F^*,v}$ for $\OO_{F^*,R}$. Thus if $R=\OO_S$ we have
\begin{equation}\label{4.6}
\fD(\OO_{F^*,S})=D(F^*)\cdot\OO_S\, .
\end{equation}

\begin{lemma}\label{L4.1}
Let $F^*=(F,\theta_F)$ be an augmented $K$-form with $F \in \OO_S[X,Y]$. Then
\begin{equation}\label{4.7}
 \fd_{ij}(\OO_{F^*,S})=[F]
\frac{[\theta_F^{(i)}-\theta_F^{(j)}]}{[1,\theta_F^{(i)}][1,\theta_F^{(j)}]}
\ \ \  (1 \leq i,j \leq r, \ i\ne j),
\end{equation}
and
\begin{equation}\label{4.8}
\prod_{1 \leq i<j \leq r} \fd_{ij}(\OO_{F^*,S})^2 =
[F]^{(r-1)(r-2)}\fD(\OO_{F^*,S}).
\end{equation}
\end{lemma}

\begin{proof}
We first prove \eqref{4.7}. Let $i,j \in \{ 1, \dots , r  \}$, $i \ne j$.
Write $F=a_0X^r+a_1X^{r-1}Y+\cdots +a_rY^r$.
Then $F=a_0\prod_{k=1}^r (X-\theta_F^{(k)}Y)$,
and so by Gauss' Lemma,
\begin{equation}\label{4.9}
[F] = [a_0]\prod_{k=1}^r [1,\theta_F^{(k)}]\, .
\end{equation}
Write
\begin{equation}\label{4.10}
\prod_{\ketsor{k=1}{k\ne i,j}}^r (X-\theta_F^{(k)}Y)=
B_0X^{r-2}+B_1X^{r-3}Y+\dots +B_{r-2}Y^{r-2}.
\end{equation}
Then $B_0=1$, and by Gauss' Lemma and \eqref{4.9},
\begin{equation}\label{4.11}
[B_0, B_1, \dots , B_{r-2}]=\prod_{k=1}^r [1,\theta_F^{(k)}]=
[F][a_0]^{-1}[1,\theta_F^{(i)}]^{-1}[1,\theta_F^{(j)}]^{-1}.
\end{equation}

Let $\{ \omega_1, \dots , \omega_r \}$ be the basis of $\OO_{F^*,S}$ given
by \eqref{1.5}. We first show that
\begin{equation}\label{4.12}
\omega_m^{(i)}-\omega_m^{(j)}=a_0B_{m-2}(\theta_F^{(i)}-\theta_F^{(j)})
\ \ \text{for} \ \ m=2, \dots , r.
\end{equation}
Write $b_k:=a_k/a_0$ for
$k = 0, \dots ,r$. Then $\prod_{k=1}^r (X-\theta_F^{(k)}Y)=
b_0X^r+b_1X^{r-1}Y+\dots +b_rY^r$.
and $a_0^{-1}\omega_m=\sum_{k=0}^{m-2} b_k\theta_F^{m-k-1}$ for $m=2\kdots r$.
Assertion \eqref{4.12} is clear for $m=2$. Let $m \geq 3$.
We have (on putting $B_{-2}=B_{-1}=0$)
$$
b_k=B_k-B_{k-1}(\theta_F^{(i)}+\theta_F^{(j)})+
B_{k-2}\theta_F^{(i)}\theta_F^{(j)} \ \ \text{for} \ \ k = 0, \dots ,r,
$$
and so
\begin{equation*}
\begin{split}
&a_0^{-1}(\omega_m^{(i)}-\omega_m^{(j)})=
\sum_{k=0}^{m-2} b_k \left( (\theta_F^{(i)})^{m-k-1}-(\theta_F^{(j)})^{m-k-1} \right) \\
&=\sum_{k=0}^{m-2} \left\{ B_k-B_{k-1}(\theta_F^{(i)}+\theta_F^{(j)})+
B_{k-2}\theta_F^{(i)}\theta_F^{(j)} \right\} \cdot
\left\{ (\theta_F^{(i)})^{m-k-1}-(\theta_F^{(j)})^{m-k-1} \right\} \\
&=\sum_{k=0}^{m-2} c_kB_k\, ,
\end{split}
\end{equation*}
where
\begin{equation*}
\begin{split}
c_{m-2}= &\theta_F^{(i)}-\theta_F^{(j)}\, ,\\
c_{m-3}= &{\theta_F^{(i)}}^2-{\theta_F^{(j)}}^2-
(\theta_F^{(i)}+\theta_F^{(j)} )(\theta_F^{(i)}-\theta_F^{(j)} )=0\, ,
\end{split}
\end{equation*}
and, if $m\geq 4$,
\begin{equation*}
\begin{split}
c_k=
{\theta_F^{(i)}}^{m-k-1}-{\theta_F^{(j)}}^{m-k-1}
-(\theta_F^{(i)}+\theta_F^{(j)}) ({\theta_F^{(i)}}^{m-k-2}-{\theta_F^{(j)}}^{m-k-2} )\qquad&\\
+\theta_F^{(i)}\theta_F^{(j)}({\theta_F^{(i)}}^{m-k-3}-{\theta_F^{(j)}}^{m-k-3}) &=0
\end{split}
\end{equation*}
for $k=0\kdots m-4$.
This implies \eqref{4.12}. By combining \eqref{4.12}, \eqref{4.11} we obtain
\begin{equation*}
\begin{split}
\fd_{ij}(\OO_{F^*,S}) &=
[\omega_2^{(i)}-\omega_2^{(j)},\dots ,\omega_r^{(i)}-\omega_r^{(j)}  ] \\
&= [a_0] \cdot [B_0, B_1, \dots , B_{r-2}]\cdot [\theta_F^{(i)}-\theta_F^{(j)}]
\\
&=[F]\frac{[\theta_F^{(i)}-\theta_F^{(j)}]}
{[1,\theta_F^{(i)}][1,\theta_F^{(j)}]}
\end{split}
\end{equation*}
which is \eqref{4.7}.

Now from
\eqref{4.6}, \eqref{1.1}, \eqref{4.9}, \eqref{4.7} we infer
\begin{equation*}
\begin{split}
\fD(\OO_{F^*,S})\OO_{L,S} 
&=[D(F)]=
[a_0^{2r-2} \prod_{1\leq i<j \leq r} (\theta_F^{(i)}-\theta_F^{(j)})^2]\\
&= [F]^{2r-2}\prod_{1\leq i<j \leq r} \left(
\frac{[\theta_F^{(i)}-\theta_F^{(j)}]}{[1,\theta_F^{(i)}][1,\theta_F^{(j)}]}
\right)^2 \\
&=[F]^{-(r-1)(r-2)}\prod_{1\leq i<j \leq r}\fd_{ij}(\OO_{F^*,S})^2,
\end{split}
\end{equation*}
which is \eqref{4.8}.
\end{proof}
\vskip 0.3cm

\begin{lemma}\label{L4.2}
Let $F^*=(F, \theta_F)$, $G^*=(G, \theta_G)$  be two augmented $K$-forms such that
\begin{align}
&F,G \in \OO_S[X,Y]; \label{4.13}\\
&\OO_{F^*,S}=\OO_{G^*,S}; \label{4.14}\\
&F^*,G^* \ \ \text{are weakly $\K$-equivalent}. \label{4.15}
\end{align}
Then $F^*,G^*$ are weakly $\OO_v$-equivalent for every $v \not\in S$.
\end{lemma}

\begin{proof}
Take $v \not\in S$. By \eqref{4.15} there are $A \in \text{GL}_2(\K)$,
$\lambda \in \K^*$ such that $G^*=\lambda F_A^*$.
Since $\OO_v$ is a principal ideal domain, there are matrices
$U_1,U_2 \in \text{GL}_2(\OO_v)$ such that
$A=U_1  \bigl( \begin{smallmatrix} a&0\\0&d \end{smallmatrix} \bigr) U_2$ with
$a,d \in \K^*$. Let $\til{F}^*:=F^*_{U_1}$, $\til{G}^*:=G^*_{U_2^{-1}}$. Then
\begin{equation}\label{4.16}
\til{F}^* \simm{\OO_v} F^*, \ \ \til{G}^* \simm{\OO_v} G^*,
\end{equation}
hence it suffices to show that $\til{F}^*, \til{G}^*$ are
weakly $\OO_v$-equivalent. Write $\til{F}^*=(\til{F}, \theta_{\til{F}})$,
$\til{G}^*=(\til{G}, \theta_{\til{G}})$. Then
$\til{G}^*=\lambda \til{F}^*_{\bigl( \begin{smallmatrix} a&0\\0&d \end{smallmatrix} \bigr)}$
which means that
\begin{equation}\label{4.17}
\til{G}(X,Y)=\lambda \til{F}(aX,dY), \ \ \ \theta_{\til{G}}=\frac{d}{a}\theta_{\til{F}}.
\end{equation}
Write $\til{F}(X,Y)=a_0X^r+a_1X^{r-1}Y+ \dots +a_rY^r$. Then $\OO_{\til{F}^*,v}$
is an $\OO_v$-module with basis
$$
\omega_1=1, \ \omega_i=\sum_{j=0}^{i-2} a_j \theta_{\til{F}}^{i-j-1} \ \ (i=2, \dots ,r).
$$
By \eqref{4.17},
$\til{G}(X,Y)=\lambda a_0 a^r X^r+\lambda a_1 a^{r-1}d X^{r-1}Y+ \dots
+\lambda a_r d^r Y^r$ and $\OO_{\til{G}^*,v}$ is an $\OO_v$-module with basis
$$
\omega_1'=1, \ \omega_i'=
\sum_{j=0}^{i-2} \lambda a_j a^{r-j} d^j
\left( \frac{d}{a}\theta_{\til{F}} \right)^{i-j-1} =
\lambda a^{r-i+1}d^{i-1} \omega_i  \ \ (i=2, \dots ,r).
$$
By \eqref{4.16}, \eqref{4.14}, \eqref{4.5} we have $\OO_{\til{F}^*,v}=\OO_{\til{G}^*,v}$.
Therefore, the matrix relating  $\{ \omega_1', \dots , \omega_r' \}$ to
$\{ \omega_1, \dots , \omega_r \}$ is in $\text{GL}_2(\OO_v)$. That is,
$$
\lambda a^{r-1}d \in \OO_v^*, \ \ \lambda a^{r-2}d^2  \in \OO_v^*, \dots ,
\lambda a d^{r-1} \in \OO_v^*,
$$
which implies $d=au$ with $u \in \OO_v^*$. Further,
$\lambda a^r=u^{-1} \lambda a^{r-1}d \in \OO_v^*$.
Inserting this into \eqref{4.17} we obtain
\begin{equation*}
\til{G}(X,Y)=\lambda \til{F}(aX,auY)=\lambda a^r \til{F}(X,uY),\ \
\theta_{\til{G}}=u\theta_{\til{F}},
\end{equation*}
which implies that $\til{F}^*,\til{G}^*$ are weakly $\OO_v$-equivalent. This proves
Lemma \ref{4.2}.
\end{proof}
\vskip 0.3cm

We now arrive at our final result:

\begin{proposition}\label{P4.3}
Let $\CC $ be a collection of augmented $K$-forms such that
\begin{align}
&F \in \OO_S[X,Y]\ \ \ \text{for every} \ \ F^*=(F,\theta_F) \in \CC ; \label{4.18}\\
&\OO_{F^*,S}=\OO_{G^*,S}\ \ \ \text{for every pair} \ \ F^*,G^* \in \CC ; \label{4.19}\\
&\text{the elements of $\CC $ are weakly $\K$-equivalent to one onother.}
\label{4.20}
\end{align}
Then if $r$ is odd, $\CC $ is contained in the union of at most $r^s$
$\OO_S$-equivalence classes, while if $r$ is even, $\CC $ is contained in
the union of at most $r^sh_2(\OO_S)$ $\OO_S$-equivalence classes.
\end{proposition}

\begin{proof}
Combine Lemmata \ref{L4.2} and \ref{L2.3}.
\end{proof}
\vskip 0.5cm

\section{Proof of Theorem \ref{T1.1}}\label{S5}

\setcounter{equation}{0}

Let $\K, S$ be as in Section \ref{S1}; thus $\#S=s$.
Let $\OO$ be an $\OO_S$-order of degree $r \geq 3$ and denote by
$K$ its quotient field.
Let $F \in \OO_S[X,Y]$ be a binary form which is irreducible in
$\K[X,Y]$ and such that $\OO_{F,S}\cong \OO$ (as $\OO_S$-algebras).
Then there is a $\theta_{F}$ such that $F(\theta_F,1)=0$, $K=\K(\theta_F)$
and such that $\omega_1, \dots , \omega_r$ given by \eqref{1.5}
form an $\OO_S$-basis of $\OO$. Thus, $F^*:=(F, \theta_F)$ is
an augmented $K$-form with $\OO_{F^*,S}=\OO$. Now it is obvious that in
order to prove Theorem \ref{T1.1} it suffices to prove the following:

\begin{proposition}\label{P5.1}
Let $\#S=s$, and  let $K$ be a finite extension of $\K$ of degree $r \geq 3$.
Let $\OO \subset K$ be an $\OO_S$-order with quotient field $K$.
Then the set of
augmented $K$-forms $F^*=(F,\theta_F)$ with
\begin{align}
&F \in \OO_S[X,Y],  \label{5.1}\\
&\OO_{F^*}=\OO \label{5.2}
\end{align}
is contained in the union of finitely many $\OO_S$-equivalence classes,
whose number is bounded above by
\begin{equation}\label{5.3}
2^{24r^3s} \ \ \text{if $r$ is odd;} \ \ \ 2^{24r^3s}h_2(\OO_S) \ \ \text{if $r$ is even}.
\end{equation}
\end{proposition}
\vskip 0.3cm

For the moment we assume $r \geq 4$. The case $r=3$ will be treated separately.
Our main tool is a result of Beukers and Schlickewei on equations in
two variables with unknowns from a multiplicative group of finite rank.
Let $\Omega$ be a field of characteristic $0$.
We endow $(\Omega^*)^2$ with coordinatewise multiplication
$(x_1,y_1) * (x_2, y_2)= (x_1x_2,y_1y_2)$; thus  $(\Omega^*)^2$ becomes
a group with unit element $(1,1)$. For $(x,y) \in (\Omega^*)^2$, $m \in \Z$
we write $(x,y)^m :=(x^m,y^m)$.

\begin{lemma}\label{L5.2}
Let $(x_1,y_1), \dots ,(x_n,y_n) \in (\Omega^*)^2$. Let
\begin{equation*}
  \begin{split}
\Gamma:=\{ (x,y)\in (\Omega^*)^2: &\exists m\in \NN, z_1, \dots ,z_n \in \Z \\
&\text{with}
\ \ (x,y)^m = (x_1,y_1)^{z_1} * \dots * (x_n,y_n)^{z_n} \}.
  \end{split}
\end{equation*}
Then the equation
\begin{equation}\label{5.4}
x+y=1 \ \ \text{in} \ \ (x,y) \in \Gamma
\end{equation}
has at most $2^{8(n+1)}$ solutions.
\end{lemma}

\begin{proof}
See \cite[Theorem 1]{BeSchl}.
\end{proof}

Let $\OO$, $K$ be as above. Choose a normal closure $L$ of $K$ over $\K$
and denote again by $\xi\mapsto \xi^{(i)}$ $(i=1\kdots r)$ the $\K$-isomorphic
embeddings of $K$ into $L$. We recall that the cross ratio
of $\al_1,\al_2,\al_3,\al_4\in\PP^1(L)$ is given by
\begin{equation}\label{5.5}
\{ \al_1,\al_2;\al_3,\al_4\}:=
\frac{(\al_1-\al_2)(\al_3-\al_4)}{(\al_1-\al_3)(\al_2-\al_4)}
\end{equation}
(with the usual adaptations if one of $\al_1\kdots\al_4$ is $\infty$
or if $\al_1\kdots\al_4$ are not all distinct).
As is well-known, cross ratios are invariant under projective 
transformations.

For an augmented $K$-form $F^*=(F,\theta_F)$
with \eqref{5.1}, \eqref{5.2}
we define the tuple of all cross ratios of 
$\theta_F^{(1)}\kdots\theta_F^{(r)}$,
\begin{equation}\label{5.6}
\Delta(F^*):=(\{\theta_F^{(i)},\theta_F^{(j)};\theta_F^{(k)},\theta_F^{(l)}  \}: \ 1 \leq i,j,k,l \leq r; \ \ i,j,k,l
\ \ \text{distinct} ) .
\end{equation}

\begin{lemma}\label{L5.3}
If $F^*$ runs through the collection of augmented $K$-forms with \eqref{5.1},
\eqref{5.2}, then $\Delta (F^*)$ runs through a collection of
cardinality at most
\begin{equation}\label{5.7}
2^{24(r^3-r^2)s}.
\end{equation}
\end{lemma}

\begin{proof}
Let $F^*=(F,\theta_F)$ be an augmented $K$-form with \eqref{5.1}, \eqref{5.2}.
Let $i,j,k,l \in \{ 1, \dots , r \}$ be distinct. We have
%$$
%(\theta_F^{(i)}-\theta_F^{(j)})(\theta_F^{(k)}-\theta_F^{(l)})+
%(\theta_F^{(j)}-\theta_F^{(k)})(\theta_F^{(i)}-\theta_F^{(l)})+
%(\theta_F^{(k)}-\theta_F^{(i)})(\theta_F^{(j)}-\theta_F^{(l)})=0
%$$
%which implies
\begin{equation}\label{5.8}
\{\theta_F^{(i)},\theta_F^{(j)};\theta_F^{(k)},\theta_F^{(l)}  \}+\{\theta_F^{(i)},\theta_F^{(l)};\theta_F^{(k)},\theta_F^{(j)}  \}=1.
\end{equation}
Write \eqref{5.8} as $x+y=1$. 
We want to apply Lemma \ref{L5.2} to \eqref{5.8} and to this end
we have to find a suitable group $\Gamma$ independent of $F^*$ such that
$(x,y)\in\Gamma$.

Fix $\theta_0$ with $\K(\theta_0)=K$. For each two-element subset $\{i,j\}$ of
$\{ 1, \dots , r \}$ define the field
\[
K^{\{i,j\}}:=\K(\theta_0^{(i)}+\theta_0^{(j)},\theta_0^{(i)}\theta_0^{(j)})\, .
\]
Thus, if $P(X,Y) \in \K[X,Y]$ is a symmetric polynomial, then
$P(\xi^{(i)},\xi^{(j)}) \in K^{\{i,j\}}$ for every $\xi\in K$.
Further, $[ K^{\{i,j\}}:\K] \leq \binom{r}{2}$.
Let $t(\{i,j\})$ denote the rank of $\OO_{K^{\{i,j\}},S}^*$, i.e., the unit group
of the integral closure of $\OO_S$ in $K^{\{i,j\}}$. Then $t(\{i,j\})$ is
equal to the number of places of $K^{\{i,j\}}$ lying above the places in $S$, minus 1.
That is,
\begin{equation}\label{5.9}
  t(\{i,j\}) \leq [ K^{\{i,j\}}:\K]s-1 \leq \binom{r}{2}s-1.
\end{equation}
There are $\ve_1^{\{i,j\}}, \dots ,\ve_{t(\{i,j\})}^{\{i,j\}} \in \OO_{K^{\{i,j\}},S}^*$
such that every element of $\OO_{K^{\{i,j\}},S}^*$ can be expressed uniquely as
\begin{equation}\label{5.10}
\zeta \prod_{m=1}^{t(\{i,j\})} (\ve_m^{\{i,j\}})^{w_m}
\end{equation}
where $\zeta \in K^{\{i,j\}}$ is a root of unity
and $w_m\in \Z$ for $m =1, \dots , t(\{i,j\})$.

Let $h$ be the least common multiple of the following integers:
the class number
of $K$; the class number of $K^{\{i,j\}}$
for each two-element subset $\{i,j\}$ of
$\{ 1, \dots , r \}$; and the number of roots of unity in $K^{\{i,j\}}$
for each two-element subset $\{i,j\}$ of $\{ 1, \dots , r \}$.

We raise the identity \eqref{4.7} to the power $2h$ to obtain something
useful. Let $i,j \in \{ 1, \dots , r \}$, $i\ne j$. First we have an
identity of fractional $\OO_{K,S}$-ideals
\begin{equation}\label{5.11}
[1, \theta_F]^{2h}=[\al_F] \ \ \text{with} \ \ \al_F \in K^*
\end{equation}
since $2h$ is a multiple of the class number of $K$.
Further, $(\theta_F^{(i)}-\theta_F^{(j)})^{2h} \in K^{\{i,j\}}$. The ideal
$\fd_{ij}(\OO_{F^*,S})^2$ is generated by elements $(\xi^{(i)}-\xi^{(j)})^2$
$(\xi \in \OO_{F^*,S})$ which belong to $K^{\{i,j\}}$. By \eqref{4.8} the
$\OO_S$-ideal $[F]$ generated by the coefficients of $F$ depends only on
$\OO_{F^*,S}$, hence by \eqref{5.2} on $\OO$. Therefore we have an identity
of fractional $\OO_{K^{\{i,j\}},S}$-ideals
\begin{equation}\label{5.12}
([F]^{-1}\fd_{ij}(\OO_{F^*,S}) )^{2h}=[\be_{ij}]\ \ \text{with} \ \ \be_{ij}
\in (K^{\{i,j\}})^*,
\end{equation}
where $\be_{ij}$ depends only on $\OO$. So in particular,
$\be_{ij}$ is independent of $F^*$. Lastly,
$\al_F^{(i)}\al_F^{(j)} \in K^{\{i,j\}}$. Now \eqref{4.7}, \eqref{5.11}, \eqref{5.12}
yield an identity of fractional $\OO_{K^{\{i,j\}},S}$-ideals
$[\theta_F^{(i)}-\theta_F^{(j)}]^{2h}=[ \al_F^{(i)}\al_F^{(j)} \be_{ij}]$, that is,
$(\theta_F^{(i)}-\theta_F^{(j)})^{2h}= \al_F^{(i)}\al_F^{(j)} \be_{ij}\eta_{ij}$
with $\eta_{ij} \in \OO_{K^{\{i,j\}},S}^*$. We can express $\eta_{ij}$ as in
\eqref{5.10}. By raising again to the power $h$, we can cancel the root of unity,
and obtain
\begin{equation}\label{5.13}
(\theta_F^{(i)}-\theta_F^{(j)})^{2h^2}=(\al_F^{(i)}\al_F^{(j)} \be_{ij})^h
\prod_{m=1}^{t(\{i,j\})} (\ve_m^{\{i,j\}})^{w_m} \ \ \text{with} \ \ w_m\in \Z.
\end{equation}
Taking any distinct $i,j,k,l \in \{ 1, \dots , r \}$,
and writing again \eqref{5.8} as $x+y=1$, it follows that
\begin{equation*}
\begin{split}
(x,y)^{2h^2}&=(\{\theta_F^{(i)},\theta_F^{(j)};\theta_F^{(k)},\theta_F^{(l)}  \},\{\theta_F^{(i)},\theta_F^{(l)};\theta_F^{(k)},\theta_F^{(j)}  \})^{2h^2} \\
& = \left(  \frac{(\theta_F^{(i)}-\theta_F^{(j)})(\theta_F^{(k)}-\theta_F^{(l)})}{(\theta_F^{(i)}-\theta_F^{(k)})(\theta_F^{(j)}-\theta_F^{(l)})} ,
\frac{(\theta_F^{(i)}-\theta_F^{(l)})(\theta_F^{(j)}-\theta_F^{(k)})}{(\theta_F^{(i)}-\theta_F^{(k)})(\theta_F^{(j)}-\theta_F^{(l)})} \right)^{2h^2} \\
& = \left(  \frac{\be_{ij} \be_{kl}}{\be_{ik} \be_{jl}} ,
\frac{\be_{il} \be_{jk}}{\be_{ik} \be_{jl}} \right)^h * (\eta_1,\eta_2)
\end{split}
\end{equation*}
where $(\eta_1,\eta_2)$ is a product of powers of
$$
(\ve_m^{\{ i,j \}} , 1 ) \ \ (1 \leq m \leq t(\{ i,j \}) ); \ \ \
(\ve_m^{\{ k,l \}} , 1 ) \ \ (1 \leq m \leq t(\{ k,l \}) ); \ \ \
$$
$$
(1, \ve_m^{\{ i,l \}}  ) \ \ (1 \leq m \leq t(\{ i,l \}) ); \ \ \
(1, \ve_m^{\{ j,k \}} ) \ \ (1 \leq m \leq t(\{ j,k \}) ); \ \ \
$$
$$
(\ve_m^{\{ i,k \}} ,\ve_m^{\{ i,k \}}  ) \ \ (1 \leq m \leq t(\{ i,k \}) ); \ \ \
(\ve_m^{\{ j,l \}} , \ve_m^{\{ j,l \}} ) \ \ (1 \leq m \leq t(\{ j,l \}) ). \ \ \
$$
It is important to notice that the terms
$\al_F^{(i)},\al_F^{(j)},\al_F^{(k)},\al_F^{(l)}$ are cancelled. Thus, in
view of \eqref{5.9}, 
$(x,y)^{2h^2}$
is a product of powers of
\begin{equation*}
\begin{split}
1+&t(\{ i,j \})+t(\{ k,l \})+t(\{ i,l \})+t(\{ j,k \})+t(\{ i,k \})+t(\{ j,l \}) \\
\leq & 1+6 \left( \binom{r}{2}s-1 \right) = 6 \binom{r}{2}s-5
\end{split}
\end{equation*}
terms which are independent of $F^*$.

Now applying Lemma \ref{L5.2} to \eqref{5.8} yields that $(x,y)$,
and so in particular $x=\{\theta_F^{(i)},\theta_F^{(j)};\theta_F^{(k)},\theta_F^{(l)}  \}$,
belongs to a set independent of $F^*$ of cardinality at most
\begin{equation}\label{5.14}
2^{8\{6 \binom{r}{2}s-5+1 \} }=2^{48 \binom{r}{2}s-32 }\, .
\end{equation}

We claim that the tuple $\Delta (F^*)$ of all cross ratios is
determined uniquely by the subtuple
\begin{equation}\label{5.15}
\til{\Delta}(F^*):=( \{\theta_F^{(1)},\theta_F^{(2)};\theta_F^{(3)},\theta_F^{(l)}  \}\  : \ l= 4, \dots ,r ).
\end{equation}
Indeed, let $\langle T \rangle$ be the unique projective transformation of $\PP^1$,
mapping $\theta_F^{(1)},\theta_F^{(2)},\theta_F^{(3)}$ to
$1,\infty,0$, respectively. Since $\langle T \rangle$ does not alter cross ratios,
for $l=4, \dots , r$ the image of $\theta_F^{(l)}$ under $\langle T \rangle$ is
$\{\theta_F^{(1)},\theta_F^{(2)};\theta_F^{(3)},\theta_F^{(l)}  \}$. 
But then it follows that  $\{\theta_F^{(i)},\theta_F^{(j)};\theta_F^{(k)},\theta_F^{(l)}  \}$ is
equal to the cross ratio of the $i$-th, $j$-th, $k$-th, $l$-th point among
$1,\infty,0,\{\theta_F^{(1)},\theta_F^{(2)};\theta_F^{(3)},\theta_F^{(4)}  \}, \dots ,\{\theta_F^{(1)},\theta_F^{(2)};\theta_F^{(3)},\theta_F^{(r)}  \}$.

So by \eqref{5.14} the total number of possibilities for
$\til{\Delta}(F^*)$, and hence that for $\Delta (F^*)$
is at most
$$
2^{\left( 48 \binom{r}{2}s-32 \right)(r-3)} \leq
2^{24(r^3-r^2)s}.
$$
This proves Lemma \ref{L5.3}.
\end{proof}
\vskip 0.3cm

\begin{lemma}\label{L5.4}
Let $F^*=(F,\theta_F)$, $G^*=(G,\theta_G)$ be two augmented $K$-forms
of degree $r \geq 3$ with \eqref{5.1}, \eqref{5.2}.

{\bf (i)} If $r=3$ then $F^*, G^*$ are weakly $\K$-equivalent.

{\bf (ii)} If $r\geq 4$ and moreover,
\begin{equation}\label{5.16}
\Delta (F^*)= \Delta (G^*),
\end{equation}
then $F^*, G^*$ are weakly $\K$-equivalent.
\end{lemma}

\begin{proof}
If $r\geq 4$ then by \eqref{5.16},
$\{ \theta_F^{(i)}, \theta_F^{(j)}; \theta_F^{(k)}, \theta_F^{(l)} \}=
\{ \theta_G^{(i)}, \theta_G^{(j)}; \theta_G^{(k)}, \theta_G^{(l)} \}$
for each distinct $i,j,k,l \in \{ 1, \dots , r \}$. This implies
that there is a
unique projective transformation $\langle T\rangle :\PP^1(L) \to \PP^1(L)$ 
with
$\langle T \rangle(\theta_F^{(i)})=\theta_G^{(i)}$ for $i=1,\dots , r$.
If $r=3$ then we simply use that there is a unique projective transformation
$\langle T \rangle:\PP^1 \to \PP^1$ defined over $\overline{\Q}$
with $\langle T \rangle(\theta_F^{(i)})=\theta_G^{(i)}$ for $i=1,2,3$.

In other words, both for $r=3$ and $r\geq 4$ there is an
up to a scalar factor unique matrix
$T = \bigl( \begin{smallmatrix} a&b\\c&d\end{smallmatrix} \bigr) \in \text{GL}_2(L)$
such that
\begin{equation}\label{5.17}
\theta_G^{(i)}=\frac{a\theta_F^{(i)}+b}{c\theta_F^{(i)}+d}
\ \ \text{for}\ \ i=1, \dots , r.
\end{equation}
We choose the first non-zero element among $a,b,c,d$ equal to 1
so that $T$ is uniquely determined.
Then for every $\tau\in\Gal(L/\K)$, the matrix
$\tau(T)=\bigl( \begin{smallmatrix} \tau(a)&\tau(b)\\\tau(c)&\tau(d)\end{smallmatrix} \bigr)$
also satisfies \eqref{5.17}
since $\tau$ permutes both sequences
$\theta_F^{(1)}, \dots ,\theta_F^{(r)}$ and
$\theta_G^{(1)}, \dots ,\theta_G^{(r)}$ in the same manner.
Hence $\tau(T)=T$ for every $\tau \in \Gal(L/\K)$
which implies $T \in \text{GL}_2(\K)$.

Write $F = a_F \prod_{i=1}^r (X-\theta_F^{(i)}Y)$,
$G = a_G \prod_{i=1}^r (X-\theta_G^{(i)}Y)$ with $a_F, a_G \in \K^*$. Thus,
\begin{equation*}
\begin{split}
G=& a_G \prod_{i=1}^r (X-\frac{a\theta_F^{(i)}+b}{c\theta_F^{(i)}+d}Y) \\
=& a_Ga_F^{-1}\left\{ \prod_{i=1}^r (c\theta_F^{(i)}+d) \right\}^{-1} F(dX-bY,-cX+aY)\\
=& a_Ga_F^{-1}\left\{ \prod_{i=1}^r (c\theta_F^{(i)}+d) \right\}^{-1}
(ad-bc)^r F_{T^{-1}}(X,Y)=\lambda F_{T^{-1}}(X,Y)
\end{split}
\end{equation*}
with $\lambda \in \K^*$, $T\in \ \text{GL}_2(\K)$, and
$\theta_G=\langle T \rangle(\theta_F)$. This implies that $F^*$, $G^*$ are
weakly $\K$-equivalent.
\end{proof}
\vskip 0.3cm

\begin{proof}[Proof of Proposition \ref{P5.1}]
Let $r \geq 3$. Put $h(r,\OO_S):=1$ if $r$ is odd,
and $h(r,\OO_S):=h_2(\OO_S)$ if $r$ is
even.
By Lemmata \ref{L5.3} and \ref{L5.4}, the collection of augmented $K$-forms
$F^*=(F,\theta_F)$ with \eqref{5.1}, \eqref{5.2} is contained in the union of
at most $2^{24(r^3-r^2)s}$ weak $\K$-equivalence classes. Together
with Proposition \ref{P4.3} this implies that the collection of
augmented $K$-forms with \eqref{5.1}, \eqref{5.2} is contained in the union of
at most
$$
2^{24(r^3-r^2)s}\cdot r^s h(r,\OO_S) \leq 2^{24r^3s}h(r,\OO_S)
$$
$\OO_S$-equivalence classes. This proves Proposition \ref{P5.1}.
\end{proof}
\vskip 0.5cm

\section{Proof of Theorem \ref{T1.2}}\label{S6}

\setcounter{equation}{0}

We keep the notation from Section \ref{S1}. Thus $\K$ is a number field and $S$ is a
finite subset of $M_{\K}$ of cardinality $s$ containing all infinite places.
Let $K$ be an extension of $\K$
of degree $r \geq 3$. Let $\fc \ne (0)$ be an ideal of $\OO_S$ and let
$S'=S \cup \{v \not\in S \ : \ |x|_v <1 \text{ for every } x \in \fc \}$.
Notice that if $F\in\F (\OO_S ,K)$ satisfies \eqref{1.10}, then
\begin{equation*}
D(F)\cdot \OO_{S'} =  \fd_{K/\K,S'}\, .
\end{equation*}
So by \eqref{1.5a}, the $\OO_{S'}$-order
associated with $F$ is $\OO_{F,S'}=\OO_{K,S'}$ (the integral closure of $\OO_{S'}$
in $K$). On applying Theorem \ref{T1.1} with $S'$ in place of $S$
and with $\OO =\OO_{K,S'}$
we infer that the set of binary forms $F\in\F (\OO_S,K)$  
with \eqref{1.10} is contained in finitely many $\OO_{S'}$-equivalence classes,
whose number is at most
\begin{equation}\label{6.3}
\begin{split}
&2^{24r^3\#S'}=2^{24r^3(s+\omega_S(\fc))} \ \ \ \ \text{if $r$ is odd,} \\
&2^{24r^3\#S'}h_2({\OO_{S'}}) \leq 2^{24r^3(s+\omega_S(\fc))}h_2({\OO_S})
\ \ \ \ \text{if $r$ is even,}
\end{split}
\end{equation}
where we have used $\#S'=s+\omega_S(\fc)$ and the obvious inequality
$h_2({\OO_{S'}}) \leq h_2({\OO_S})$. 

In particular, the binary forms $F\in\F (\OO_S,K)$ with
\eqref{1.10} lie in finitely many $\K$-equivalence classes,
whose number is bounded above by \eqref{6.3}.
By multiplying this quantity with the upper bound \eqref{3.100}
from Proposition \eqref{P3.7}, (ii) we obtain an upper bound for
the number of $\OO_S$-equivalence classes of binary forms under
consideration which is precisely the upper bound from Theorem \ref{T1.2}.
This completes our proof.\qed
\vskip0.5cm

\section{Proof of Theorem \ref{T1.3}}\label{S7}

\setcounter{equation}{0}

To prove Theorem \ref{T1.3}, we need a further Proposition
on resultant equations which can be
regarded as a quantitative version of Lemma 1 of Evertse and Gy\H ory
\cite{EGy8}.

For the moment, let $K_0, K_1$ be two (not necessarily distinct) extensions of $\K$
of degrees $r_0, r_1$, respectively, such that $r_0 \geq 3$.
Let $L$ be a normal closure over $\K$ of the compositum of $K_0,K_1$.
Below, by $[a_1\kdots a_m]$ we will denote the fractional $\OO_{L,S}$-ideal
generated by $a_1\kdots a_m$, and by $[f]$ the fractional $\OO_{L,S}$-ideal
generated by the coefficients of a given polynomial $f$.

Using the notation of Theorems \ref{T1.2} and \ref{T1.3}, fix a binary form
$F_0 \in \F(\OO_S,K_0)$,
and consider the binary forms $F_1 \in \F(\OO_S, K_1)$.

\begin{proposition}\label{P7.1}
Up to multiplication by $S$-units, there are at most
$$
2^{24r_0r_1s}
$$
binary forms $F_1 \in \F(\OO_S, K_1)$ which satisfy
\begin{equation}\label{7.1}
R(F_0,F_1) \in \OO_S^*.
\end{equation}
\end{proposition}

\begin{proof}
Take $F_1\in\F(\OO_S ,K_1)$ with \eqref{7.1}.
By assumption, for $i=0,1$ we have that $F_i \in \OO_S[X,Y]$,
$F_i$ is irreducible over $\K$,
and there is a $\theta_i$ satisfying $F(\theta_i, 1)=0$ and $\K(\theta_i)=K_i$.
We can write
$$
F_i(X,Y)=a_i\prod_{k=1}^{r_i} (X-\theta_i^{(k)}Y)\quad (i=0,1),
$$
where $a_i\in\K^*$, and where
$\theta_i^{(1)}\kdots \theta_i^{(r_i)}$ are
the conjugates of $\theta_i$ in $L$, for $i=0,1$.
By Gauss' Lemma we have
\begin{equation}
\label{7.3}
[1] \supseteq [F_i]=[a_i]\prod_{k=1}^{r_i} [1,\theta_i^{(k)}]\quad (i=0,1).
\end{equation}
Using \eqref{7.1} and expression \eqref{3.3a} for the resultant, we get
\begin{equation*}
\begin{split}
[1] &=[R(F_0,F_1)]=
[a_0]^{r_1}[a_1]^{r_0}
\prod_{k=1}^{r_0}\prod_{l=1}^{r_1}[\theta_0^{(k)}-\theta_1^{(l)}]
\\
&\subseteq
\prod_{k=1}^{r_0}\prod_{l=1}^{r_1}
\frac{[\theta_0^{(k)}-\theta_1^{(l)}]}{[1,\theta_0^{(k)}][1,\theta_1^{(l)}]}\, .
\end{split}
\end{equation*}
In combination  with the obvious inclusions
$[\theta_0^{(k)}-\theta_1^{(l)}]\subseteq [1,\theta_0^{(k)}][1,\theta_1^{(l)}]$
this gives
\begin{equation}\label{7.6}
[\theta_0^{(k)}-\theta_1^{(l)}]=[1,\theta_0^{(k)}][1,\theta_1^{(l)}]
\ \ \mbox{for $k=1\kdots r_0$, $l=1\kdots r_1$.}
\end{equation}
Meanwhile, we have shown also that the inclusions in \eqref{7.3}
are equalities, i.e.,
\begin{equation}
\label{7.4}
[F_i] =[a_i]\prod_{k=1}^{r_i} [1,\theta_i^{(k)}]=[1]\quad (i=0,1).
\end{equation}

We proceed similarly as in the proof of Lemma \ref{L5.3}.
Define the fields $K_{i1}:=\K(\theta_0^{(i)},\theta_1)=K_0^{(i)}K_1$
$(i=1\kdots r_0)$.
Denote by $h$ the least common multiple of the class numbers of
$K_0$, $K_1$, $K_{11}\kdots K_{r_0,1}$
and of the numbers of roots of unity of
$K_{11}\kdots K_{r_0,1}$.
By our choice of $h$,
there are $\al_0 \in K_0^*$ such that $[1, \theta_0]^h = [\al_0]$, and
$\al_1 \in K_1^*$ such that $[1, \theta_1]^h = [\al_1]$.
Then by \eqref{7.6},
$$
[\theta_0^{(i)}-\theta_1]^h=[\al_0^{(i)}][\al_1] \ \
\text{for} \ \ i=1, \dots , r_0,
$$
that is
$$
(\theta_0^{(i)}-\theta_1)^h=\al_0^{(i)} \al_1 \eta_i,
$$
where $\eta_i\in\OO_{K_{i1},S}^*$ (i.e., the unit group of the integral closure
of $\OO_S$ in $K_{i1}$).
Let $\ve_{i1}, \dots , \ve_{is_i}$ be a system of fundamental units
of $\OO_{K_{i1},S}^*$. Then $\eta_i$ is a product of a root of unity in
$K_{i1}$ and of powers of $\ve_{i1}, \dots , \ve_{is_i}$ and so, by our
choice of $h$,
$$
(\theta_0^{(i)}-\theta_1)^{h^2}=(\al_0^{(i)})^h \al_1^h
\ve_{i1}^{w_{i1}} \dots \ve_{is_i}^{w_{is_i}}
$$
with $w_{i1}, \dots , w_{is_i} \in \Z$.

Pick distinct subscripts $i, j, k \in \{ 1, \dots , r_0\}$ and
consider the identity
$$
\frac{(\theta_0^{(i)}-\theta_0^{(j)})}{(\theta_0^{(i)}-\theta_0^{(k)})}\cdot
\frac{(\theta_0^{(k)}-\theta_1)}{(\theta_0^{(j)}-\theta_1)}
+
\frac{(\theta_0^{(j)}-\theta_0^{(k)})}{(\theta_0^{(i)}-\theta_0^{(k)})}\cdot
\frac{(\theta_0^{(i)}-\theta_1)}{(\theta_0^{(j)}-\theta_1)}
=1.
$$
This can be written as $x+y=1$, where
$$
(x,y)^{h^2}=(a,b)* \prod_{q=1}^{s_k} (\ve_{kq},1)^{w_{kq}}
* \prod_{q=1}^{s_i} (1,\ve_{iq})^{w_{iq}}
* \prod_{q=1}^{s_j}(\varepsilon_{jq},\varepsilon_{jq})^{-w_{jq}}
$$
with
$$
(a,b)=\left( \Big(
\frac{(\theta_0^{(i)}-\theta_0^{(j)})}{(\theta_0^{(i)}-\theta_0^{(k)})}
\Big)^{h^2} \Big(\frac{\alpha_{0}^{(k)} }{\alpha_{0}^{(j)} }\Big)^h\, ,\,
\Big(\frac{(\theta_0^{(j)}-\theta_0^{(k)})}{(\theta_0^{(i)}-\theta_0^{(k)})}
\Big)^{h^2}\Big(\frac{\alpha_{0}^{(i)} }{\alpha_{0}^{(j)} }\Big)^h\right)\, .
$$
Notice that
$$
s_i\leq ([K_{i1}:\K]s)-1\leq r_0r_1s-1
$$
and similarly for $s_j$ and $s_k$.
So, by Lemma \ref{5.2} the number of possibilities for $(x,y)$ is at most
$$
2^{8(s_i+s_j+s_k+1)+8}\leq 2^{24r_0r_1s}\, .
$$
This gives at most $2^{24r_0r_1s}$ possibilities for $\theta_1$.
But, by \eqref{7.4}, the ideal $[a_1]$
is uniquely determined once $\theta_1$ is uniquely determined
and moreover, $a_1\in\K^*$.
So $a_1$ is uniquely determined up to a factor from $\OO_S^*$.
We infer that up to multiplication
by some factor from $\OO_S^*$, for $F_1$ there are at most $2^{24r_0r_1s}$
possibilities.
\end{proof}
\vskip 0.3cm

\begin{proof}[Proof of Theorem \ref{T1.3}]
Let $K_0, K_1, \dots , K_t$ be (not necessarily distinct) extensions of $\K$
of degrees $r_0, r_1, \dots , r_t$, respectively, such that $r_0 \geq 3$.
Let $F\in \F(\OO_S, K_0, \dots , K_t)$ be a binary form
with the property \eqref{1.12}. There are binary forms $F_0\kdots F_t$
with $F=F_0\cdots F_t$ and with $F_i\in\F (\OO_S,K_i)$ for $i=0\kdots t$.
So in particular, $F_i\in\OO_S[X,Y]$ for $i=0\kdots t$.
Let $S'$ denote the union of $S$ and the places $v\not\in S$
such that $|x|_v<1$ for every $x\in\fc$.
Then
\begin{equation*}
D(F)\cdot \OO_{S'}=\fd_{K_0/\K ,S'}\cdots \fd_{K_t/\K ,S'}\, .
\end{equation*}
Now by expressing $D(F)$ as in \eqref{3.4a}, and using
$R(F_i,F_j)\in\OO_{S'}$ $(0\leq i<j\leq t)$
and the inclusions
\begin{equation*}
D(F_i)\cdot\OO_{S'}\subseteq \fd_{K_i/\K ,S'}\quad (i=0\kdots t)
\end{equation*}
(which follow from (ii) of Lemma \ref{L3.1}),
we obtain
\begin{eqnarray}
\label{7.7}
&&D(F_0)\cdot \OO_{S'}=\fd_{K_0/\K, S'}\, ,
\\
\label{7.7a}
&&R(F_0,F_i)\in\OO_{S'}^*\quad (i=0\kdots t)\, .
\end{eqnarray}

We apply now Theorem \ref{T1.2} to \eqref{7.7} with $S$ replaced by $S'$;
we obtain that $F_0$
is contained in the union of at most
\begin{equation}\label{7.8}
2^{24r_0^3(\# S')}h(r_0, \OO_{S'})
\leq 2^{24r_0^3(s+\omega_S(\fc))}h(r_0, \OO_S)
\end{equation}
$\OO_{S'}$-equivalence classes. Here we have used that
$\# S' = s+\omega_S(\fc))$ and
$h(r_0, \OO_{S'})\leq h(r_0, \OO_S)$.

Fix one of these $\OO_{S'}$-equivalence classes, and pick 
from this class a representative
$F_0 \in \F (\OO_S ,K_0)$ with \eqref{7.7}.
Consider all tuples $(F_1\kdots F_t)$ of binary forms
with $F_i\in\F (\OO_S ,K_i)$ for $i=1\kdots t$ and with \eqref{7.7a}.
Proposition \ref{P7.1} gives that for given $F_0$
there are, up to $S'$-unit factors,
at most
$$
2^{24r_0(r_1+\cdots +r_t)(s+\omega_S(\fc))}
$$
such tuples $(F_1\kdots F_t)$.

Combining this with the upper bound \eqref{7.8} for the number of
$\OO_{S'}$-equiv\-al\-ence classes of binary forms $F_0\in\F (\OO_S ,K_0)$
with \eqref{7.7}, we infer that up
to $\OO_{S'}$-equivalence, and up to
an $\OO_{S'}$-unit factor,
there are at most
\begin{equation}\label{7.8a}
\begin{split}
&2^{24r_0^3(s+\omega_S(\fc))}h(r_0,S)
\cdot2^{24r_0(r_1+\cdots +r_t)(s+\omega_S(\fc ))}
\\
&\qquad\quad =
2^{24r_0(r_0^2+r_1+\cdots +r_t)(s+\omega_S(\fc ))}h(r_0,\OO_S)
\end{split}
\end{equation}
binary forms $F=F_0\cdots F_t\in\F (\OO_S,K_0\kdots K_t)$
with \eqref{1.12}.
That is, there are binary forms $G_1\kdots G_m\in\F (\OO_S,K_0\kdots K_t)$,
with $m$ bounded above by the quantity in \eqref{7.8a},
such that every binary form $F\in\F (\OO_S,K_0\kdots K_t)$
with \eqref{1.12} is $\OO_{S'}$-equivalent to $\ve G_i$
for some $i\in\{ 1\kdots m\}$ and $\ve\in\OO_{S'}^*$.
But $\ve$ can be
written in the form $\ve_1^{w_1}\cdots\ve_{s'}^{w_{s'}}\eta^r$, where
$s'=\# S'=s+\omega_S(\fc )$, 
$\ve_1\kdots\ve_{s'}$ are generators of $\OO_{S'}^*$, 
$w_1\kdots w_{s'}\in\{ 0\kdots r-1\}$ and $\eta\in\OO_{S'}^*$.
Since $G_i$ is $\OO_{S'}$-equivalent to
$\eta^rG_i=(G_i)_{\bigl(\begin{smallmatrix}\eta&0\\0&\eta\end{smallmatrix}}\bigr)$,
we have in fact that every binary form $F$ under consideration is
$\OO_{S'}$-equivalent to $\ve_1^{w_1}\cdots\ve_{s'}^{w_{s'}}G_i$, 
with $w_1\kdots w_{s'}\in\{ 0\kdots r-1\}$
and with $i\in\{ 1\kdots m\}$.
Assuming as we may in view of Theorem \ref{T1.2}
that $r_1+\cdots +r_t\geq 1$,
it follows that the binary forms
$F\in\F (\OO_S, K_0\kdots K_t)$ with \eqref{1.12}
lie in at most
\begin{equation*}
\begin{split}
&\Big( r\cdot 2^{24r_0(r_0^2+r_1+\cdots +r_t)}\Big)^{(s+\omega_S(\fc ))}h(r_0,\OO_S)
\\
&\qquad
\leq \Big( r\cdot 2^{24(r-1)((r-1)^2+1)}\Big)^{(s+\omega_S(\fc ))}h(r_0,\OO_S)
\end{split}
\end{equation*}
and so in at most
\begin{equation}\label{7.9}
2^{24r^3(s+\omega_S(\fc ))}h(r_0,\OO_S)
\end{equation}
$\OO_{S'}$-equivalence classes.

By (ii) of Proposition \ref{P3.7}, 
the binary forms $F\in\F (\OO_S ,K_0\kdots K_t)$ with \eqref{1.12}
lie in finitely many $\OO_S$-equivalence classes whose product is bounded
above by the product of \eqref{7.9} and of \eqref{3.100}.
Since this is precisely the bound of Theorem \ref{T1.3},
this completes our proof.
\end{proof}
\vskip 0.5cm

\section{Lower bounds}\label{S8}

\setcounter{equation}{0}

We present some examples, showing that the results mentioned in Section \ref{S1}
are in certain respects close to best possible. 
\vskip 0.2cm

First let $K$ be a finite extension of $\K$ of even degree $r \geq 4$.
Let $S$ be a finite subset of $M_{\K}$ such that $S$ contains all infinite
places.
We show that there are infinitely many $\OO_S$-orders
$\OO$ with quotient field $K$, such that the collection of augmented $K$-forms
$F^*=(F,\theta_F)$ with $F\in \OO_S[X,Y]$ and $\OO_{F^*,S}=\OO$ cannot be
contained in fewer than $h_2(\OO_S)$ $\OO_S$-equivalence classes. 
Since
each binary form $F\in\F (\OO_S,K)$ gives rise to at most
$r$ augmented $K$-forms $F^*=(F,\theta_F)$, it follows 
that
the set of forms $F \in \F (\OO_S ,K)$ with
$\OO_{F,S}\cong\OO$ cannot be
contained in fewer than $r^{-1}h_2(\OO_S)$ $\OO_S$-equivalence classes.
This shows that the factor $h_2(\OO_S)$ in the upper bound of Theorem \ref{T1.1}
is necessary.

Pick any augmented $K$-form $F^*=(F,\theta_F)$ with $F\in \OO_S[X,Y]$. Let
$\A$ be any ideal of $\OO_S$ such that $\A^2$ is principal.
The ideal $\A$ can be generated by two elements,
$\A=[\al, \be]$, say.
Let $\A^2=[\lambda]$. Then there are $\xi, \eta \in \OO_S$ such that
$\xi\al^2-\eta\be^2=\lambda$.
Define
$$
F^*_{\A}:=\lambda^{-r/2}
F^*_{\bigl( \begin{smallmatrix} \al&\be \\ \eta\be&\xi\al\end{smallmatrix} \bigr)}\, .
$$

We first show that $F^*_{\A}=(F_{\A},\theta_{F_{\A}})$ with
$F_{\A} \in \OO_S[X,Y]$, and $\OO_{F_{\A}^*,S}=\OO_{F^*,S}$. 
Pick $v \not\in S$.
Then there is $\mu \in \OO_v$ such that in $\OO_v$ we have the identity of
ideals $[\al, \be]=[\mu]$. We now get
$$
F_{\A}=\lambda^{-r/2} F(\al X+\be Y, \eta \be X+\xi \al Y)=
\lambda^{-r/2} \mu^r
F(\frac{\al}{\mu} X+\frac{\be}{\mu} Y, \frac{\eta \be}{\mu} X+\frac{\xi \al}{\mu} Y).
$$
Since $[\mu^2]=[\lambda]$ in $\OO_v$ we have $\lambda^{-r/2} \mu^r \in \OO_v^*$.
Further,
$$
\det \begin{pmatrix}\frac{\al}{\mu}& \frac{\be}{\mu}\\
\frac{\eta \be}{\mu}&\frac{\xi \al}{\mu}  \end{pmatrix}=
\frac{\xi\al^2-\eta\be^2}{\mu^2}=\frac{\lambda}{\mu^2} \in \OO_v^*.
$$
Hence $F^*_{\A}$, $F^*$ are weakly $\OO_v$-equivalent. This implies
$F_{\A} \in \OO_v[X,Y]$. Further by \eqref{4.5},
$\OO_{F_{\A}^*,v}=\OO_{F^*,v}$
where $\OO_{F_{\A}^*,v}, \OO_{F^*,v}$ are the localizations at $v$ of
$\OO_{F_{\A}^*,S}, \OO_{F^*,S}$. This holds for every $v \not\in S$.
Hence $F_{\A} \in \OO_S[X,Y]$ and $\OO_{F_{\A}^*,S}= \OO_{F^*,S}$.

We now show that if $\A_1$, $\A_2$ are two ideals of $\OO_S$ such that
$\A_1^2$, $\A_2^2$ are principal and $\A_1$, $\A_2$ do not belong to
the same ideal class, then the augmented $K$-forms $F_{\A_1}^*, F_{\A_2}^*$
constructed above are not $\OO_S$-equivalent. Thus, the 
collection of augmented
$K$-forms $F_{\A}^*$ such that $\A$ is an ideal of $\OO_S$ for which
$\A^2$ is principal cannot be contained in fewer than $h_2(\OO_S)$
$\OO_S$-equivalence classes.

For $i=1,2$ let $\A_i=[\al_i, \be_i]$ be an ideal of $\OO_S$, suppose that
$\A_i^2=[\lambda_i]$ is principal, and choose $\xi_i, \eta_i \in \OO_S$
such that $\xi_i\al_i^2-\eta_i\be_i^2=\lambda_i$ for $i=1,2$.
Define $F_{\A_i}^*:=
\lambda_i^{-r/2}
F^*_{\bigl( \begin{smallmatrix} \al_i&\be_i \\ \eta_i\be_i&\xi_i\al_i\end{smallmatrix} \bigr)}$
$(i=1,2)$.
Suppose that $F_{\A_2}^*=(F_{\A_1}^*)_U$ for some $U \in \text{GL}_2(\OO_S)$.
Then by (ii) of Lemma \ref{L2.1}, there is $\rho \in \K^*$ such that
$\bigl( \begin{smallmatrix} \al_2&\be_2 \\ \eta_2\be_2&\xi_2\al_2\end{smallmatrix} \bigr)=
\rho \bigl( \begin{smallmatrix} \al_1&\be_1 \\ \eta_1\be_1&\xi_1\al_1\end{smallmatrix} \bigr)U$, and
$\rho^r=(\lambda_1\lambda_2^{-1})^{r/2}$. Hence $[\rho]^r=(\A_1\A_2^{-1})^r$
which implies  $\A_1=\rho \A_2$. So $\A_1, \A_2$ lie in the same ideal
class. This proves our assertion.
\vskip 0.2cm

Now let $(K_0, \dots , K_t)$ be a sequence of finite extensions of $\K$ such that
$\sum_{i=0}^t [K_i:\K]=:r \geq 3$. We show that there are infinitely many
ideals $\fc$ of $\OO_S$ such that the collection of binary forms
$\F (\OO_S, K_0\kdots K_t)$ with \eqref{1.12}
cannot be contained in fewer than
$C\times N_S(\fc)^{2/r(r-1)}$ $\OO_S$-equivalence classes,
where $C$ is some positive constant.

Fix $\til{F}\in\F (\OO_S,K_0\kdots K_t)$ with $D(\til{F})\ne 0$.
Extend this to
an augmented $(K_0, \dots , K_t)$-form
$\til{F}^*=(\til{F}, \theta_{0,\til{F}}, \dots , \theta_{t,\til{F}})$.
Let $a\in \OO_S$, $a\ne 0$.
For $\be \in \OO_S$ define
$$
\til{F}^*_{\be}:=\til{F}^*_{\bigl( \begin{smallmatrix} 1&\be \\ 0&a \end{smallmatrix} \bigr)}=
(\til{F}_{\be}, \theta_{0,\til{F}_{\be}}, \dots , \theta_{t,\til{F}_{\be}})
\ \ \text{with} \ \ \til{F}_{\be}=\til{F}(X+\be Y, aY).
$$
Now if $\be_1, \be_2 \in \OO_S$ are such that $\til{F}^*_{\be_1}, \til{F}^*_{\be_2}$
are $\OO_S$-equivalent, then
$\til{F}^*_{\bigl( \begin{smallmatrix} 1&\be_1 \\ 0&a \end{smallmatrix} \bigr)}=
\til{F}^*_{\bigl( \begin{smallmatrix} 1&\be_2 \\ 0&a \end{smallmatrix} \bigr)U}$
for some matrix
$U \in \text{GL}_2(\OO_S)$. According to Lemma \ref{L2.1}, (ii), this implies
$\bigl( \begin{smallmatrix} 1&\be_1 \\ 0&a \end{smallmatrix} \bigr)^{-1}
\bigl( \begin{smallmatrix} 1&\be_2 \\ 0&a \end{smallmatrix} \bigr) \in \text{GL}_2(\OO_S)$
and therefore,
$(\be_1-\be_2)/a \in \OO_S$.

Consequently, the
augmented $(K_0, \dots , K_t)$-forms $\til{F}^*_{\be}$ $(\be \in \OO_S)$
cannot be contained in the union of fewer than
$\#\OO_S/[a]=N_S(a)$ $\OO_S$-equivalence classes.

Notice that
$\til{F}_{\beta}\in \F (\OO_S, K_0\kdots K_t)$ for $\beta\in\OO_S$.
By (ii) of Lemma \ref{L3.1}, there is an ideal $\fc_0$ of $\OO_S$
such that
$[D(\til{F})]=  \fc_0^2 \fd_{K_0/\K,S} \dots \fd_{K_t/\K,S}$.
Put $\fc :=a^{\frac{1}{2}r(r-1)}\fc_0$. Then by \eqref{1.3},
$\til{F}_{\be}$ satisfies \eqref{1.12} with this $\fc$.

Since there 
at most $r^{t+1}$ different augmented forms $\til{F}^*_{\be}$
coming from the same binary form $\til{F}_{\beta}$,
it follows that for each ideal $\fc$ as constructed above,
the set of binary forms $F\in\F (\OO_S ,K_0\kdots K_t)$
with \eqref{1.12} cannot be contained in the union of fewer than
$$
r^{-t-1}N_S(a)=r^{-t-1}N_S(\fc_0)^{\frac{-2}{r(r-1)}}N_S(\fc)^{\frac{2}{r(r-1)}}=:
C\times N_S(\fc)^{\frac{2}{r(r-1)}}
$$
$\OO_S$-equivalence classes.

\end{document}